\documentclass{comp12pt}

\usepackage[english]{babel}

\usepackage{isolatin1}

\usepackage{amsthm}

\usepackage{amsmath}

\usepackage{amssymb} 

\usepackage{amscd}

\usepackage[all]{xy}

\usepackage{times}

\usepackage{mathrsfs}

\renewcommand{\theenumi}{(\roman{enumi})}

\renewcommand{\thesection}{\arabic{section}}

\renewcommand{\theequation}{\arabic{equation}}

\swapnumbers

\theoremstyle{plain}
\newtheorem{theorem}[subsubsection]{Theorem}
\newtheorem{proposition}[subsubsection]{Proposition}
\newtheorem{lemma}[subsubsection]{Lemma}
\newtheorem{corollary}[subsubsection]{Corollary}

\newtheorem{definition}[subsubsection]{Definition}

\theoremstyle{definition}
\newtheorem{remark}[subsubsection]{Remark}

\newtheorem{void}[subsubsection]{}

\numberwithin{equation}{subsubsection}

\makeatletter
\renewcommand{\subsubsection}
{\@startsection {subsubsection}{3}{\z@ }{-3.25ex\@plus -1ex \@minus -.2ex} 
{-1.5ex}{\normalfont \normalsize \bfseries }}
\makeatother

\renewcommand{\mathcal}{\mathscr}

\renewcommand{\Bbb}{\mathbb}

\newcommand{\goth}{\mathfrak}

\renewcommand{\bf}{\mathbf}
\renewcommand{\cal}{\mathcal}

\DeclareMathOperator{\End}{End} 

\DeclareMathOperator{\GL}{GL} 
\renewcommand{\H}{\mathrm{H}} 
\DeclareMathOperator{\Hom}{Hom} 
\DeclareMathOperator{\SheafHom}{\mathbf{Hom}}

\DeclareMathOperator{\Pic}{Pic}

\newcommand{\red}{\mathrm{red}} 
\newcommand{\Sch}{\mathbf{Sch}}
\DeclareMathOperator{\Spec}{Spec}
\DeclareMathOperator{\hgt}{ht}  
\DeclareMathOperator{\rk}{rk}   

\newcommand{\isomto}{\overset{\sim}{\rightarrow}}
\newcommand{\isomfrom}{\overset{\sim}{\leftarrow}}
\newcommand{\surjto}{\twoheadrightarrow}
\newcommand{\injto}{\hookrightarrow}

\newcommand{\CC}{\Bbb{C}}
\newcommand{\FF}{\Bbb{F}} 
\newcommand{\GG}{\Bbb{G}} 
\newcommand{\QQ}{\Bbb{Q}}
\newcommand{\RR}{\Bbb{R}}
\newcommand{\ZZ}{\Bbb{Z}}

\newcommand{\kk}{\Bbbk} 

\newcommand{\A}{\bf{A}}
\newcommand{\K}{\bf{K}}

\newcommand{\R}{\bf{R}}
\newcommand{\C}{\bf{C}}

\renewcommand{\k}{\overline{k}} 

\newcommand{\E}{\bf{E}} 
\newcommand{\F}{\bf{F}}

\renewcommand{\L}{\mathcal{L}} 
\newcommand{\M}{\mathcal{M}}

\renewcommand{\O}{\mathcal{O}} 

\newcommand{\n}{\goth{n}} 
\newcommand{\m}{\goth{m}} 
\newcommand{\p}{\goth{p}} 
\newcommand{\q}{\goth{q}} 

\newcommand{\Y}{\underline{Y}} 
\newcommand{\Yo}{\underline{Y_0}} 
\newcommand{\Yi}{\underline{Y_1}} 

\renewcommand{\o}{\overline}

\begin{document}

\title
 {Drinfeld modular curves have many points}

\author{Lenny Taelman}
\email{lenny@math.rug.nl}
\address{IWI\\RuG Groningen\\Postbus 800\\9700AV Groningen\\Nederland}

          
\thanks{The author would like to thank Jan Van Geel.}

\begin{abstract}
 An algebraic (smooth, projective) curve over a finite field of
 $q$ elements can have at most $(\sqrt{q}-1+o(1))g$ rational
 points for a growing genus $g$. If $q$ is a square, then
 this bound is sharp. Many examples have been found of 
 families of curves of increasing genus attaining this bound
 (by, amongst others, M.~A.~Tsfasman, S.~G.~Vladut,
 T.~Zink, Y.~Ihara, A.~Garcia, H.~Stichtenoth, N.~Elkies).
 These are typically towers of reductions of modular curves
 or Drinfeld modular curves. In the former case,
 it has been proven by Y.~Ihara that all towers of 
 reductions of modular curves of a certain type attain this
 bound. In the latter case, this has been proven only
 for particular examples, even though, as remarked by
 Noam Elkies, a similar argument should work.
 
 This text applies the argument of Ihara to the case
 of Drinfeld modular curves. The bulk of work in doing
 this is working out the theory of (reductions of)
 Drinfeld modular curve, as a lot of properties of
 these curves that are needed for Ihara's argument
 to work are absent from the literature. As a consequence,
 the text is almost completely self-contained.
 
 The main result can be stated as follows.
 Let $X$ be a smooth projective and absolutely irreducible
 curve over the finite field $\FF_q$ of $q$ elements. Let
 $\infty$ be a point of $X$ and denote by $\A$ the ring of
 functions on $X$ that are regular outside $\infty$. 
 Let $\n$ be a principle ideal in $\A$ and $\p$ a prime
 ideal not containing $\n$. Let $q^m$ be the cardinality of
 the residue field of $\p$. Denote by $X_0(\n)$ the 
 compactification of the Drinfeld modular variety over
 $\A$ that classifies the Drinfeld $\A$-modules of rank $2$ 
 together with a submodule isomorphic to $\A/\n\A$.
 Then the reduction modulo $\p$ of $X_0(\n)$ is
 a curve of genus $g$ having at least
  $ (q^m-1)(g-1) $ rational points over the quadratic extension
 of the residue field of $\p$.
\end{abstract}

\maketitle

\tableofcontents

\section{Notation}

This thesis is written in the language of schemes and representable
functors, so let us fix some notation and make some conventions. 

All schemes are assumed to be locally noetherian. We do not distinguish between
rings or algebras (both taken to be commutative and with unit element,
unless otherwise stated) and affine schemes. Often we do not write
$\Spec(R)$ but just $R$ for the spec of the ring $R$. It should be
clear from the context (usually from the direction of the arrows) if we are
speaking ``algebra'' or ``scheme''. 

Categories are usually not abbreviated. For example: $(\text{abelian
groups})$. Except for the category $(\Sch/S)$ of schemes over a base scheme
$S$ or $(\Sch)$ for $(\Sch/\ZZ)$. 

Isomorphisms are denoted by $\isomto$, $\isomfrom$ and $\cong$. Canonical
isomorphisms simply by $=$.

\section{Line bundles in positive characteristic}

\subsection{Review}
In this chapter we study additive group schemes in positive characteristic.
All schemes in this chapter will be schemes over a finite field $\FF_q$
of characteristic $p$. We will study additive groups equipped with not only
an $\FF_p$-linear structure, which they carry in a natural way, but also an
$\FF_q$-linear one. Let us first recall some 
facts about line bundles and invertible sheaves. We fix a base scheme
$S\to\FF_q$.

\begin{void}
A line bundle $G$ over $S$ is a group scheme which is zariski locally
isomorphic to the additive group $\GG_a$ on $S$. We call $G$ together with
an action of $\FF_q$, that is a ring morphism
 $$ \FF_q \to \End(G) $$
into the ring of group scheme endomorphisms, an \emph{$\FF_q$-line
bundle}. We will from now on assume that all line bundles are equipped with
such an $\FF_q$-structure.

This means that $G$ is not only a commutative group object but also an
$\FF_q$-vector space object in the category $(\Sch/\FF_q)$.
A morphism $G \to H$ of $\FF_q$-line bundles over $S$ is said to be
$\FF_q$-linear (or short: an $\FF_q$-morphism) if it commutes with the action
of the finite field $\FF_q$. This is the same as saying that $G \to H$ is a
morphism of functors, where $G$ and $H$ are interpreted as contravariant
functors from $(\Sch/S)$ to $(\text{vector spaces}/\FF_q)$.

A trivialization provides an open affine covering of $S$ with spectra of rings
$R_i$, such that $G\times_S R_i \cong R_i[X]$, where the addition on
$R_i[X]$ is given by
 $$ R_i[X] \to R_i[X]\otimes_{R_i}R_i[X]=R_i[Y,Z]  : X \mapsto Y+Z $$
and the $\FF_q$-linear structure by
 $$ \FF_q \to \End(R_i[X]) : \lambda \mapsto (X \mapsto \lambda X). $$
\end{void}

\begin{void}
An \emph{invertible sheaf} on $S$ is a coherent $\O_S$-module sheaf
$\L$, locally free of rank 1. It is called invertible because the dual sheaf
  $$ \L^{-1} := \SheafHom(\L,\O_S) $$
satisfies the following canonical isomorphism:
  $$ \L \otimes \L^{-1} = \O_S $$
\end{void}

\begin{void}
Invertible sheaves correspond to line bundles and
vice versa. There is a one-one correspondence between isomorphism classes
of line bundles and isomorphism classes of invertible sheaves. This
can be established in two ways.

Given an invertible sheaf $\L$ on $S$, the contravariant functor 
  $$ (f:T\to S) \mapsto \H^0( T , f^\ast(\L)  ) $$
of $(\Sch/S)$ to $(\text{commutative groups})$ is
representable by a line bundle $G$, we call it $G=\GG_{a,\L}$.

Alternatively, we can construct the scheme $G=\GG_{a,\L}$ as the affine
scheme on $S$, corresponding to the symmetric algebra of the sheaf
$\L^{-1}$ (\cite{EGA1} I 9.4). 

The converse construction is obtained by restricting the functor of points
of $G/S$ to open immersions $U\to S$. The restricted functor is representable
by an invertible sheaf $\L$.
\end{void}

\subsection{Homomorphisms}

This section is largely based on section 1.2 of Thomas Lehmkuhl's
habilitationsschrift \cite{Lehmkuhl00}. From here on, we suppress the $\FF_q$ and
just talk about ``line bundles'' whenever we mean ``$\FF_q$-line
bundles''.

\begin{void}
A line bundle on $S$ can be trivialized by covering $S$ with open affines.
On such an affine, corresponding to a ring $R$, the restriction of the line
bundle is isomorphic to $\GG_{a,R} = R[X]$, with the usual additive group
structure on $R[X]$. On this group scheme we can easily find a ring of
$\FF_q$-endomorphisms, namely by adjoining the $q$-th power frobenius
$\tau$ to the ring $R$. $\tau$ does not commute with multiplication by an
element of $R$. The resulting ring, of which the elements are endomorphisms
working on the left, is the skew polynomial ring $R\{\tau\}$ with
commutation rule $ r^q \tau = \tau r $. 
\end{void}

\begin{proposition}
  $R\{\tau\}$ is the full ring $\End_{\FF_q}(\GG_{a,R})$ of $\FF_q$-linear
  endomorphisms of $\GG_{a,R}$.
\end{proposition}

\begin{proof}
  An $\FF_q$-linear morphism $ \GG_{a,R}\to\GG_{a,R} $ corresponds to a
  ring morphism $ R[X] \gets R[X] $, sending $X$ to an $\FF_q$-linear, additive
  polynomial $ f(X) = a_0 + a_1 X + \cdots + a_n X^n \in R[X] $. Given such
  an $\FF_q$-linear $f$, we want to prove that for all indices $i$ that are
  not a power of $q$ the coefficient $a_i$ vanishes. This is a purely
  algebraic exercise. For details, see for example (\cite{Lehmkuhl00} 1.2.1).
\end{proof}

\begin{void}
We now use this affine result to determine the abelian group
$\Hom_{\FF_q}(\GG_{a,\L},\GG_{a,\M})$ of $\FF_q$-morphisms of commutative
group schemes. Consider the abelian group of global sections
  $$ \H^0(S,\bigoplus_{n\geq0} \L\otimes\M^{-q^i}). $$
This group is canonically isomorphic to
  $$ \bigoplus_{n\geq0}^{loc. fin.} \H^0( S , \L\otimes\M^{-q^n}) $$ 
of which the elements are formal power series in $\tau$
  $$ a_0 + a_1\tau + a_2\tau^2 + \cdots $$
with $ a_i \in \H^0( S, \M\otimes\L^{-q^i}) $ such that locally only a
finite number of terms are non-zero.
A term $a_i\tau^i$ of an element of such a power series
defines an $\FF_q$-linear morphism of invertible sheaves by
composing the $i$-th iterated $q$-th power frobenius
  $$ \tau^i: \L \to \L^{q^i} : s \mapsto s \otimes s \cdots \otimes s $$
with multiplication by a global section:
  $$ a_i: \L^{q^i} \to \M : t \mapsto t \otimes a_i \in 
     \L^{q^i}\otimes\L^{-q^i}\otimes\M = \M. $$
Because we demand our power series to be locally finite, we can locally sum
these morphisms to $\FF_q$-linear morphisms, and glue them to form a global
linear morphism of line bundles.
If $\L = \M$, this group of endomorphisms is actually a ring of
endomorphisms, multiplication is given by
  $$ a_i\tau^i b_j\tau^j = a_i \otimes b_j^{q^i} \tau^{i+j} $$
Now we claim that, as in the local case, this group is the full
group of $\FF_q$-morphisms.
\end{void}

\begin{proposition}
  Let $\L$ and $\M$ be invertible sheaves on the $\FF_q$-scheme
  $S$. We have a canonical isomorphism of groups
  $$ \Hom_{\FF_q}(\GG_{a,\L},\GG_{a,\M}) = \
     \bigoplus_{n\geq0}^{loc. fin.} \H^0( S , \L\otimes\M^{-q^n}) $$
  and a canonical isomorphism of rings
  $$ \End_{\FF_q}(\GG_{a,\L}) = \
     \bigoplus_{n\geq0}^{loc. fin.} \H^0( S , \L^{1-q^n}) $$
\end{proposition} 

\begin{proof}
  Given an $\FF_q$-morphism $\GG_{a,\L}\to\GG_{a,\M}$, we can write it
  locally (on the affines of a trivializing covering) in the desired form.
  Since this is done uniquely, it follows from the uniqueness on the
  intersections that we can glue the coefficients to form global sections
  over $S$ of the sheaves $\L\otimes\M^{-q^n}$. For the converse, we have
  already seen that every series of sections of the desired form induces
  canonically an $\FF_q$-morphism of line bundles. 
\end{proof}

We are especially interested in finite morphisms. We restrict our attention 
to connected base schemes, to make notions as rank defined globally.

\begin{proposition}\label{finitemorphisms}
  Given a non-trivial $\FF_q$-morphism $\xi:\GG_{a,\L}\to\GG_{a,\M}$ of
  line bundles over the connected scheme $S$, the following are equivalent:
  \begin{enumerate}
    \item $\xi$ is finite
    \item $\xi$ is finite and flat, of rank $q^n$, for some $n$
    \item there exists an $n$ such that
          $\xi = x_0 + x_1\tau + x_2\tau^2 + \cdots$ with $x_n$ invertible,
          i.e., nowhere vanishing, and all the $x_i$ with $i>n$ are
	  nilpotent, i.e., everywhere vanishing.
  \end{enumerate}
  If $\L=\M$, and if the above conditions are satisfied, then there is a unique
  $\FF_q$-automorphism $\sigma=\sum{s_i\tau^i}$ of $\GG_{a,\L}$ with
  $s_0=1$
  such that
    $$ \sigma^{-1}\circ\xi\circ\sigma = x_0+z_1\tau+\cdots+z_n\tau^n $$
  where $z_n$ is invertible.
\end{proposition}

\begin{proof}
  (i) $\Rightarrow$ (iii). If $\xi$ is finite, it is automatically flat
  since flatness is local and every finite $R$-morphism $R[X]\to R[X]$
  defines on the target the structure of a free module over the source.
  Because $\xi$ is flat it has a locally constant rank. If we write $\xi$
  in the usual form
    $$ \xi = x_0 + x_1\tau + x_2\tau^2 + \cdots $$
  this means that in every point of $S$, some $x_n$ is non-zero, and all
  higher coefficients are zero. Moreover, since $S$ is connected, this $n$
  is a constant.
  
  (iii) $\Rightarrow$ (ii). True on all open affine subschemes of $S$,
  hence also globally true.
  
  (ii) $\Rightarrow$ (i). A fortiori.
  
  For the last claim, assume again $S$ is a ring $R$. Write $\xi$ as a
  finite polynomial in $\tau$
    $$ \xi=x_0+x_1\tau+x_2\tau^2+\cdots+x_n\tau^n +\cdots+x_m\tau^m $$ 
  with all coefficients $x_i$ for $n < i \leq m$ generating a nilpotent
  ideal $\goth{n}$ of $R$. Let us assume $\n^{2^s}=0$.
  Remark that
    $$ \left(1+\frac{x_m}{x_n^{q^{m-n}}}\tau^{m-n}\right)^{-1} \equiv 
       \left(1-\frac{x_m}{x_n^{q^{m-n}}}\tau^{m-n}\right)
       \mod \n^{2^{s-1}} $$
  And that the coefficients in $\tau^i$ for $i\geq m$ of
    $$ \left(1+\frac{x_m}{x_n^{q^{m-n}}}\tau^{m-n}\right)^{-1}
       \xi 
       \left(1+\frac{x_m}{x_n^{q^{m-n}}}\tau^{m-n}\right) $$
  vanish modulo $\n^{2^{s-1}}$. So by induction on the degree of $\xi$
  modulo $\n^{2^{s-1}}$, we find an automorphism $\sigma$ such that
    $$ \xi' = \sigma^{-1}\xi\sigma = \sum_{0}^{m'}{x'_i\tau^i} $$
  has coefficients $x'_i$ for $n < i \leq m'$ generating a nilpotent ideal
  $\n'$ satisfying $\n'^{2^{s-1}}=0$. Now by decreasing induction on $s$ we
  find the existence of a $\sigma$ satisfying the desired property.
  
  For the uniqueness, it is sufficient to show that $\sigma$ must be
  trivial if $\xi$ has already the desired form with zero coefficients
  for $i>n$. Because $\sigma$ must be invertible, it looks like
    $$ \sigma = 1 + s_1\tau + s_2\tau^2 + \cdots + s_t\tau^t $$
  where all the $s_i$ are nilpotent. Let $\m$ be the nilpotent ideal
  generated by them, and let $s$ be minimal such that $\m^{2^s}=0$. We
  will apply the same trick as for the existence. If $s=0$, we are done. 
  Now assume $s>0$ and let $i$ be the highest 
  index such that $s_i$ is non-zero modulo $\m^{2^{s-1}}$. By evaluating
  the terms in $\tau^{n+i}$ in
    $$ (x_0+\cdots+x_n\tau^n)\sigma \equiv \sigma(z_0+\cdots+z_n\tau^n)
       \mod \m^{2^{s-1}} $$
  we find that $s_i$ vanishes modulo $\m^{2^{s-1}}$, a contradiction.

  Now for a general $S$, consider an affine covering. On every affine
  we find a $\sigma$, and these glue because of the uniqueness on the
  intersections of the affines. Moreover, the global $\sigma$ is unique
  because it is unique on affines.
\end{proof}

\begin{corollary}\label{corfin}
  A non-trivial finite $\FF_q$-morphism of line bundles over $S$ is an
  epimorphism of $S$-schemes.
\end{corollary}

\begin{proof}
  It follows from the above explicit description that it is surjective
  (because it is surjective on geometric points), and hence faithfully
  flat. We also see that it is quasi-compact. These two properties imply
  the result since faithfully flat and quasi-compact morphisms are descent
  morphisms, and therefore strict epimorphisms, hence a fortiori
  epimorphisms.
\end{proof}

\begin{void}
Given a finite $\FF_q$-morphism $\xi:\GG_{a,\L}\to\GG_{a,\M}$ of line
bundles over the connected scheme $S$, we define its \emph{height} to be
the smallest index $h$ for which the coefficient $x_h$ is non-zero.
The \emph{derivative at $0$} is defined to be the
coefficient $x_0$. We
use the notations $\rk$, $\hgt$ and $\partial$ for rank, height and
derivative at $0$ respectively:
  $$ \rk(x_0 + \cdots + x_n\tau^n + \text{ nilpotents }) = q^n $$
  $$ \hgt(x_h\tau^h + x_{h+1}\tau^{h+1} + \cdots) = h $$
  $$ \partial(x_0 + x_1\tau + \cdots) = x_0. $$
The operators $\rk$ and $\hgt$ satisfy the following relations:
  \begin{eqnarray*}
    && \rk(\xi \circ \alpha) = \rk(\xi) \rk(\alpha) \\
    && \rk(\alpha + \beta) \leq \max(\rk(\alpha),\rk(\beta))\\
    && \hgt(\xi \circ \alpha) = \hgt(\xi) + \hgt(\alpha)\\
    && \hgt(\alpha + \beta) \geq \min(\hgt(\alpha),\hgt(\beta)).
  \end{eqnarray*}
\end{void}

\subsection{Quotients by finite flat subgroup schemes}

In this section, we fix a base scheme $S/\FF_q$. 

\begin{void}
Let us first recall the definition of a quotient group scheme
(see for example \cite{Tate97}). We say that a short sequence
  $$ 0 \to X \to Y \to Z \to 0 $$
of flat commutative group schemes over $S$, where $X$ is finite is
\emph{exact} if for every group scheme $T/S$, the induced sequences
  $$ 0 \to \Hom(T,X) \to \Hom(T,Y) \to \Hom(T,Z) $$
and
  $$ 0 \to \Hom(Z,T) \to \Hom(Y,T) \to \Hom(X,T) $$
are exact sequences of (abstract) commutative groups. If this is the
case, we also say that $Z$ is the quotient of $Y$ by $X$. Remark that
$Z$ does not necessarily represent the functor $ T \mapsto Y(T)/X(T) $,
but $Y(T)/X(T)$ does map injectively (and functorially) into $ Z(T) $.
\end{void}

We call a finite flat subgroup scheme $G \to \GG_{a,\L}$ a finite flat
$\FF_q$-subgroup scheme if it has an $\FF_q$-linear structure, that is
inherited from the one on $\GG_{a,\L}$ by base change. Above a connected
component of $S$, $G$ has a constant rank $r$, which is always a power of $p$,
and even of $q$ if $G$ is an $\FF_q$-subgroup scheme. Clearly the kernel of
an $\FF_q$-morphism of line bundles over $S$ is a finite flat
$\FF_q$-subgroup scheme, but the converse is also true.

\begin{theorem}\label{quotientbundle}
  Given a finite flat $\FF_q$-subgroup scheme $G$ of $\GG_{a,\L}$, then
  there is, unique up to a unique isomorphism, an $\FF_q$-morphism of line
  bundles
    $$ \xi: \GG_{a,\L} \to \GG_{a,\M} $$
  such that the sequence
    $$ 0 \to G \to \GG_{a,\L} \to \GG_{a,\M} \to 0 $$
  is an exact sequence of group schemes over $S$. Moreover, if $S$ is
  connected and the rank of
  $G$ is denoted by $r$, then $\M$ is isomorphic to $\L^{\otimes r}$ and the
  rank of $\xi$ is also $r$.
\end{theorem}

\begin{proof}
  Without loss of generality, we may assume $S$ is connected. We know that
  the quotient $H=\GG_{a,\L}/G$ exists, as a group scheme. Let
  us see what it looks like locally on $S$. 
  
  On an open affine $R$ of $S$, over which $\GG_{a,\L}$ is trivial, we have
   $$\GG_{a,\L}\times_S R \cong R[X]$$
  Moreover, since $G$ is flat and $\FF_q$-linear, there is a monic
  $\FF_q$-linear polynomial $P(X)=\sum{a_iX^{q^i}}$ such that 
   $$ G\times_S R \cong R[X]/(P(X)). $$
  We see that above $R$, the morphism $\GG_{a,\L}\to H$ is given by
   $$ R[X] \leftarrow R[X] : P(X) \leftarrowtail X $$
  Since the rank of $G$ is the degree of $P(X)$, the claim follows from
  this local discussion.
\end{proof}

\section{Drinfeld modules}

Fix a nonsingular projective curve with finite field of constants $\FF_q$ of
characteristic $p$. Choose a closed point $\infty$ on it and let $\A$ be the
coordinate ring of the affine curve obtained by removing $\infty$. Then
$\A$ is a dedekind ring but need not be a principal ideal domain.
We write $\K$ for
its quotient field, that is the function field of the curve. Let
$\R$ be the completion of $\K$ at the place $\infty$, and $\C$ the completion
of the algebraic closure of $\R$. Then $\C$ will be algebraically closed
as well. For a prime ideal $\p$ of $\A$, we write $\A \to \kk(\p)$
for the corresponding residue field and $\A \to \overline{\kk(\p)}$ for a fixed
algebraic closure. For every non-zero ideal $\m$ of $\A$, $\A/\m$ is a
finite-dimensional vector space over $\FF_q$. We call its dimension the
degree of $\m$.

There is a well-known deep analogy of $\A$, $\K$, $\R$ and $\C$ with,
respectively, $\ZZ$, $\QQ$, $\RR$ and $\CC$. As usual, we will call the former
the ``function field setting'', and the latter the ``classical setting''.
Drinfeld modules
(called elliptic modules by Drinfel'd in his paper \cite{Drinfeld74E})
are analogs over function fields of abelian varieties over number fields. 
We will follow \cite{Deligne69} in calling drinfeld modules of rank 2, roughly
corresponding to elliptic curves, elliptic modules.

\subsection{Drinfeld modules over $\A$-schemes}

Let $\gamma:S\to\A$ be an $\A$-scheme. Roughly, a drinfeld module over $S$
is a nontrivial algebraic action of $\A$ on a line bundle over $S$.

\begin{definition}
  A \emph{drinfeld module of rank $r>0$ over $S$} is a pair $\E=(\L,\phi)$
  of an invertible sheaf $\L$ and an $\FF_q$-linear ring morphism
    $$ \phi: \A\to\End_{\FF_q}(\GG_{a,\L}) : a\mapsto\phi_a $$
  such that for every non-zero $a\in\A$ of degree $d$ we have:
  \begin{enumerate}
    \item $\partial \phi_a = \gamma^\sharp(a) $
    \item $\phi_a$ is finite of rank $q^{rd}$.
  \end{enumerate}
  A morphism between two drinfeld modules of rank $r$ over $S$ is a morphism
  of group schemes that is compatible with the actions of $\A$.
\end{definition}

\begin{void}
If $T \to S$ is an $S$-scheme, then $\GG_{a,\L}(T)$ inherits functorially
the $A$-action of $\E$. We write $\E(T)$ for the (abstract) $\A$-module of
$T$-valued points. Thus, we are using the symbol $\E$ both for the
drinfeld module itself and for the functor
  $$ \E : (\Sch/S) \to (\A\text{-modules}) $$
it represents. 

If we fix a scheme $f: T \to S$, and restrict the functor $\E$ to the
category of $T$-schemes, it is represented by a new drinfeld module $\E_T/T$,
which is obtained from $\E/S$ by base change
  $$ \E_T = \left(f^\ast\L,\psi\right) $$
where
  $$ \psi : a\mapsto \psi_a = \phi_a \times_S T. $$
\end{void}

We have already seen how a single finite endomorphism of a line bundle can
be normalized, now we will do the same with all endomorphisms $\phi_a$
simultaneously. The resulting object is a standard drinfeld module.

\begin{definition}
  A \emph{standard drinfeld module of rank $r$ over $\gamma:S\to\A$} is a pair
  $\E=(\L, \phi)$ of an invertible sheaf $\L$ and an $\FF_q$-linear
  ring morphism
    $$ \phi: \A \to \End_{\FF_q}(\GG_{a,\L}) : a \mapsto \phi_a $$
  of the form
    $$ \phi_a = \gamma^\sharp(a) + a_1\tau + \cdots + a_{rd}\tau^{rd} $$
  for non-zero $a\in\A$, where $d$ is the degree of $a$ and where $a_{rd}$
  does not vanish on $S$.
\end{definition}

Clearly a standard drinfeld module is a drinfeld module. The
converse is also true in the following sense:

\begin{proposition}
  Every drinfeld module of rank $r$ over $S$ is isomorphic to a
  standard drinfeld module of the same rank over $S$. Moreover, this
  isomorphism is unique op to a unique isomorphism.
\end{proposition}

\begin{proof}
  Choose a non-zero $a\in\A$, say of degree $d$. Then there is, by
  \ref{finitemorphisms}, a unique isomorphism $\E\to\F=(\L,\psi)$ such that
  $\psi_a$ is in standard form. Now assume for some non-zero $b\in\A$ the
  endomorphism $\psi_b$ of $\GG_{a,\L}$ is not standard. This implies there
  exists an affine subscheme of $S$, corresponding to a ring $R$, such that
    $$ \psi_b\otimes R = \beta_0 + \ldots + \beta_i\tau^i $$
  is not standard, so $\beta_i$ is nilpotent. If we write down the
  restriction of $\psi_a$ to $R$:
    $$ \psi_a\otimes R = \alpha_0 + \ldots + \alpha_{rd}\tau^{rd} $$
  where $\alpha_{rd}$ is a unit in $R$, then we read from
  $\psi_a\psi_b=\psi_b\psi_a$ that
    $$ \beta_i^{q^{rd}} = \beta_i \alpha_{rd}^{q^i-1} $$
  This contradicts that $\beta_i$ is nilpotent.
\end{proof}

\begin{void}
Let $\E=(\L,\phi)$ be a drinfeld module over the $\A$-scheme $S$.
This scheme can be covered by open affines that trivialize the line bundle
$\L$. So we can see a drinfeld module over a scheme as a glued
collection of drinfeld modules over affine schemes with trivial invertible
sheaf. In standard form, such a drinfeld module of rank $r$ over the
$\A$-algebra $ \gamma^\sharp: \A \to R $ is just a ring morphism
$ \psi: \A \to \End(\GG_{a,R}) $ that maps a non-zero $a$ to an endomorphism
$$ \psi_a = \gamma^\sharp(a) + a_1\tau + \cdots + a_{rd}\tau^{rd} $$
where $d$ is the degree of $a$, the sections $a_n$ are elements of $R$,
the leading coefficient $a_{rd}$ is a unit and $\tau$ is just the
$q$-th power frobenius endomorphism of $\GG_{a,R}$.
\end{void}

\begin{void}
We can also see a drinfeld module over $S\to\A$ as an algebraic collection
of drinfeld modules over the residue fields of $S$. We will often give more
precise results for base schemes $S$ which are fields, or even
algebraically closed fields.
\end{void}

\subsection{Isogenies}

\begin{void}
Let again  $\gamma:S\to\A$ be an $\A$-scheme, and let $\E=(\L,\phi)$ and
$\F=(\M,\psi)$ be drinfeld modules of rank $r$ over $S$. An \emph{isogeny}
from $\E$ to $\F$ is a finite morphism of group schemes
$\xi:\GG_{a,\L} \to \GG_{a,\M}$ compatible with the $\A$-action on both
objects. This means for every $a\in\A$ the following diagram commutes:
\[\begin{CD}
  \GG_{a,\L}   @>\xi>>    \GG_{a,\M} \\
  @V\phi_aVV              @VV\psi_aV \\
  \GG_{a,\L}   @>\xi>>    \GG_{a,\M} 
\end{CD}\]
\end{void}

\begin{void}
A special isogeny of $\E$ onto itself is \emph{multiplication by $a$} for a
non-zero $a \in \A$. We denote it by $[a]=\phi_a$. The ring of all
isogenies of $\E \to \E$ is denoted by $\End(\E)$. As in the classical
setting, we say that $\E$ has complex multiplication if the map 
$\A\to\End(\E)$ is not surjective. The kernel of an isogeny $\xi:\E\to\F$
is a finite flat group scheme over $S$. We denote the kernel of $\xi$ by
$\E[\xi]$, the kernel of $[a]$ by $\E[a]$. The intersection of the
kernels of $[a]$, that is their fiber product over $\GG_{a,\L}$,
for all non-zero $a$ in some non-zero ideal $\n$ of $\A$ is denoted by
$\E[\n]$. We call $\E[\n]$ the \emph{$\n$-torsion scheme of $\E$}.
\end{void}

\begin{proposition}\label{torsionscheme}
  Let $\n$ be a non-zero prime ideal of $\A$ with prime factorization
  $$ \n = \q_1^{n_1} \q_2^{n_2} \cdots \q_m^{n_m}.$$
  \begin{enumerate}
    \item $\E[\n]$ is a finite flat group scheme over $S$. It is \'{e}tale
          over $S$ if $\n$ is disjoint from the characteristic of $S$.
    \item $\E[\n]$ has rank $|\A/\n|^r$ over $S$.
    \item $\phi$ induces an action of $\A$ on $\E[\n]$ which factors through
          $\A/\n$.
    \item There is a canonical isomorphism 
          $$ \E[\n] = \E[\q_1^{n_1}] \times_S \E[\q_2^{n_2}] \\ 
            \times_S \cdots \times_S \E[\q_m^{n_m}]. $$
  \end{enumerate}
\end{proposition}

\begin{proof}
  For all the claims it is sufficient to prove them over affine open
  subschemes. 

  Let $\E/R$ be a drinfeld module of rank $r$ over the algebra $\A\to R$
  with trivial line bundle $\GG_{a,R} = R[X]$. Assume first that $\n$ is a
  principal ideal, generated by $a\in \A$ of degree $\deg(a)=d$. We have
   $$ \E[\n] = R[X]/(\gamma(a)X + a_1X^q + \ldots + a_{rd}X^{q^{rd}}) $$
  This is finite and flat because $a_{rd}$ is invertible in $R$. By the
  Jacobi criterion, it is \'{e}tale whenever $\gamma(a)$ is invertible.
  This is precisely the case when $\n$ is disjoint from the characteristic.
  From the
  above explicit discription we read that the rank of $\E[\n]$ is $|\A/\n|^r$.
  The action of $\phi$ factors through $\A/\n$ because $\n$ is in the kernel
  of it. 
  Because fiber products of finite flat or finite \'{e}tale schemes are
  again finite flat, respectively finite \'{e}tale, the proposition is also
  true for non-principal ideals $\n$.
  The last claim follows from the third one and the chinese remainder
  theorem.
\end{proof}

We would also want to know when a finite flat sub-module is the kernel of an
isogeny: when does the quotient of a drinfeld module by such a sub-module
exist and is itself a drinfeld module? Lemma \ref{quotientbundle} tells us
that the quotient line bundle exists, so the remaining question is if it
can be equipped with a drinfeld $\A$-action.

\begin{proposition}\label{quotientmodule}
  Let $\E=(\L,\phi)$ be a drinfeld module of rank $r$ over the
  connected  $\A$-scheme $\gamma:S\to \A$. Let $G$ be a closed
  $\A$-invariant subgroup scheme of $\GG_{a,\L}$, finite and flat over $S$.
  There exists a drinfeld module $\F$ over $S$, of equal rank, and an
  isogeny $\xi:\E\to\F$ with kernel $G$ if and only if, for all $a\in\A$
  the condition
  $$ \gamma^\sharp(a^{q^h}) = \gamma^\sharp(a) $$ 
  holds, where $h$ is the height of $G$.
\end{proposition}

\begin{proof}
  The line bundle of $\F$ must be the line bundle $\M$ of theorem
  \ref{quotientbundle}. We now try to construct a drinfeld module
  $\F=(\M,\psi)$ with the desired properties. Because the quotient map
  $\xi:\GG_{a,\L}\to\GG_{a,\M}$ is an epimorphism (\ref{corfin}), and
  because $\xi$ has to be a morphism of drinfeld modules, there can only be
  one choice for $\psi$. Namely, for all $a \in \A$, $\psi_a$ is determined by
  $$ \xi \circ \phi_a = \psi_a \circ \xi $$
  Because $G$ is $\A$-invariant, a $\psi_a$ satisfying this relation exists,
  and $\psi$ is automatically a ring morphism of $\A$ into
  $\End(\GG_{a,\M})$. 

  First we consider an affine scheme $S$ with trivial line bundle $\L$.
  This is just the additive group $\GG_{a,R}$ of an $\A$-algebra
  $\gamma^\sharp: \A \to R$. Theorem \ref{quotientbundle} tells us that
  the line bundle $\M$ will also be trivial. Therefore the quotient
  $\xi:\GG_{a,\L}\to\GG_{a,\M}$ is an endomorphism of
  $\GG_{a,R}$, which has the form
  $\xi = x_h\tau^n + x_{h+1}\tau^{h+1} + \cdots$ with non-zero $x_h$. 
  Now, by comparing coefficients in $ \xi\circ\phi_a = \psi_a\circ\xi $,
  we find a necessary and sufficient condition for $\F=(\M,\psi)$ to
  be a drinfeld module of rank $r$, namely:
  $$ \forall a \in \A: \gamma^\sharp(a^{q^h}) = \gamma^\sharp(a). $$
  
  Now consider a general $\A$-scheme $S$. The only condition to be checked
  is the compatibility of the induced action on the quotient with the
  structure morphism $S\to \A$, this can be done locally.
\end{proof}

Let $\A\to k$ be an $\A$-field of characteristic $\p$, and $\A\to\k$ an
algebraic closure of $k$. For such base schemes we can say a bit more
about torsion and isogenies.

\begin{theorem}
  If $\n=\q_1^{n_1} \q_2^{n_2} \cdots \q_m^{n_m}$ is the
  factorization of $\n$ in prime ideals, then the $\n$-torsion scheme
  of the rank $r$ drinfeld module $\E/k$ is, as $\A$-module
  $$ \E[\n](k) =
     \E[\q_1^{n_1}](k) \oplus \cdots \oplus \E[\q_m^{n_m}](k) $$
  Now let $\q$ be a non-zero prime ideal of $\A$, then either
  \begin{enumerate}
    \item $\q \neq \p$ and $\E[\q^n](\k) \cong (\A/\q^n)^r$ for all $n > 0$
    \item $\q = \p$ and $\E[\q^n](\k) \cong (\A/\q^n)^h$ for all $n > 0$
      where $0 \leq h \leq r-1$ is a constant.
  \end{enumerate}
\end{theorem}

In the second case, we call an elliptic module (that is a drinfeld module
of rank 2) \emph{ordinary} if $h=1$ and \emph{super-singular} if $h=0$.

\begin{proof}
  The first statement is just \ref{torsionscheme}. Now let $m$ be the
  order of the class of $\q$ in $\Pic(\A)$, this means $m$ is the smallest
  positive integer such that $\q^m$ is a principal ideal, generated by,
  say, $a$. For positive integers $l$ we have
  $$ \phi_a \E[\q^{lm}](\k) = \E[\q^{lm-m}](\k). $$
  But because of \ref{torsionscheme} we also have
  $$ |\E[\q^{lm}](\k)| = |\E[a^l](\k)| = q^{ls} $$
  where $s = \deg(\phi_a)$ if $\q \neq \p$ and
  $s = \deg(\phi_a) - \hgt(\phi_a)$ if $\q = \p$. From this already
  the theorem follows for $n$ divisible by $m$. Now, for general $n$ use
  the filtration
  $$ \E[\q^{lm}](\k) \hookrightarrow
     \E[\q^{lm+1}](\k) \cdots \hookrightarrow
     \E[\q^{lm+m}](\k) $$
  and the action of an element in $\q$, not in $\q^2$.
\end{proof}

\section{Level structures}

Drinfel'd introduced level structures for drinfeld modules in his paper
\cite{Drinfeld74E}. Using the same idea in the classical context, a good definition
of level structures for ``bad primes'' was then finally given. A very good
book about the classical theory of moduli problems, which focuses on the bad
primes using Drinfel'd's definition is ``Arithmetic moduli of elliptic
curves'' (\cite{Katz85}).

\subsection{$\Gamma(\n)$-structures}

\begin{void}
Let $\E/S/\A=(\L,\phi)$ be an elliptic module over the $\A$-scheme $S$.
Consider a morphism of $\A$-modules
  $$ \alpha: (\A/\n)^2 \to \GG_{a,\L}(S) $$
Every element $x$ of $(\A/\n)$ defines an $S$-rational point $\alpha(x)$,
and therefore also a cartier divisor $[\alpha(x)]$ on $\GG_{a,\L}$. If the
sum of all these divisors is the cartier divisor corresponding to the
closed subscheme $\E[\n]$ of $\GG_{a,\L}$, we say $\alpha$ is a
\emph{$\Gamma(\n)$-structure} on $\E$. 

Shorter: a $\Gamma(\n)$-structure on $\E$ is an $\A$-morphism
  $$ \alpha: (\A/\n)^2 \to \GG_{a,\L}(S) $$
inducing an equality of cartier divisors
  $$ \sum_{x \in (\A/\n)^2}{[\alpha(x)]} = \E[\n] $$
on $\GG_{a,\L}$.

$\Gamma(\n)$-structures commute with base change. An elliptic module over $S$
only allows $\Gamma(\n)$-structures if its $\n$-torsion points are
$S$-rational, therefore every elliptic module can be equipped with
$\Gamma(\n)$-structures after finite faithfully flat 
base change. If $\n$ is disjoint from the characteristic of $S$,
$\Gamma(\n)$-structures exist even after finite \'{e}tale base change.

Let $\n=\q_1^{n_1}\q_2^{n_2}\cdots\q_m^{n_m}$ be the factorization of $\n$
into prime ideals. Then giving a $\Gamma(\n)$-structure is the same as giving
a $\Gamma(\q_i^{n_i})$-structure for every $i$.
\end{void}

\begin{void}
Now assume our base scheme is an algebraically closed $\A$-field $\A\to
\overline{k}$. If $\n$ is coprime to the characteristic of
$\overline{k}$, then a $\Gamma(\n)$-structure on $\E/\overline{k}$ is just
an isomorphism of $\A$-modules
  $$ (\A/\n)^2 \isomto \E[\n](\overline{k}) \cong (\A/\n)^2 $$ 
If $\n=\q^n$, where the prime $\q$ is the characteristic of $\overline{k}$,
then a $\Gamma(\n)$-structure is a surjection of $\A$-modules
  $$ (\A/\n)^2 \surjto \E[\n](\overline{k}) \cong
     \left\{ \begin{array}{ll}
       \A/\n &(\text{ordinary }\E) \\
       0     &(\text{super-singular }\E) \end{array} \right. $$
\end{void}

\subsection{$\Gamma_1(\n)$-structures}

\begin{void}
A \emph{$\Gamma_1(\n)$-structure} on $\E/S$ is a morphism of $\A$-modules
  $$ \alpha: \A/\n \to \GG_{a,\L}(S) $$
such that the effective cartier divisor
  $$ \sum_{x \in \A/\n}{[\alpha(x)]} $$
is a subgroup scheme of $\E[\n]$. 

$\Gamma_1(\n)$-structures commute with base change. An elliptic module that
allows $\Gamma(\n)$-structures also allows $\Gamma_1(\n)$-structures, but
not necessarily vice versa. 

Just as with $\Gamma(\n)$-structures, we can factor
$\Gamma_1(\n)$-structures: giving a $\Gamma_1(\n)$-structure is the same as
giving a $\Gamma_1(\q_i^{n_i})$-structure for every $i$, where
$\prod_i{\q_i^{n_i}}$ is the prime factorization of $\n$.
\end{void}

\begin{void}
Let $\E/\overline{k}$ be an elliptic module over an algebraically closed
field $\A\to\overline{k}$. If $\n$ is disjoint from the characteristic of
$\overline{k}$, then a $\Gamma_1(\n)$-structure is an injection of
$\A$-modules:
  $$ \A/\n \injto \E[\n](\overline{k}) \cong (\A/\n)^2. $$
If $\n=\q^n$, where the prime $\q$ is the characteristic of $\overline{k}$,
then a $\Gamma_1(\n)$-structure is either an isomorphism
  $$ \A/\n \isomto \E[\n](\overline{k}) $$
or the zero morphism
  $$ \A/\n \to 0 \injto \E[\n](\overline{k}) $$
where the former can occur only if $\E$ is ordinary, and the latter occurs for
both ordinary and super-singular elliptic modules.
\end{void}

\subsection{$\Gamma_0(\n)$-structures}

The definition of a $\Gamma_0(\n)$-structure is a bit more technical
because we have to consider $S$-rational subgroups of which the elements
themselves need not be $S$-rational.

\begin{void}
A \emph{$\Gamma_0(\n)$-structure} on $\E/S$ is a finite flat subgroup
scheme
  $$ H \subset \E[\n] $$
with an induced action of $\A/\n$, of constant rank $\#(\A/\n)$ over $S$, and 
\emph{cyclic}. This last condition means there is a finite, faithfully flat
base-change $T\to S$, and a point $P \in (\GG_{a,\L}\times_S T)(T)$
such that we have on $\GG_{a,\L}\times_S T$ an equality of cartier divisors
  $$ H_T = \sum_{x\in\A/\n}{[xP]}. $$ 
We call such a point $P$ a \emph{base point}.
 
$\Gamma_0(\n)$-structures are stable under base change because finite,
faithfully flat morphisms are. An elliptic module $\E/S$ need not allow
$\Gamma_0(\n)$-structures. If $\E/S$ has a $\Gamma_1(\n)$-structure
$\alpha$, it has an induced $\Gamma_0(\n)$-structure by taking $H$ to be
the effective cartier divisor
  $$ \sum_{x\in\A/\n}{[\alpha(x)]}. $$

Factor $\n$ as $\prod_i{\q_i^{n_i}}$. Given a $\Gamma_0(\n)$-structure $H$, we
can construct $\Gamma_0(\q_i^{n_i})$-structures $ H_i $ by taking the
$\q^i$-torsion in $H$:
  $$ H_i = H[\q_i^{n_i}] = \bigcap_{a\in\q_i^{n_i}}{\ker([a]:H\to H)}. $$
Conversely, given $\Gamma_0(\q_i^{n_i})$-structures $ H_i $ for all $i$, we
construct a $\Gamma_0(\n)$-structure $H$ as
  $$ H = H_1 \times_S H_2 \times_S \cdots
     \subset \E[\q_1^{n_1}] \times_S \E[\q_2^{n_2}] \times_S \cdots
     = \E[\n]
  $$
These constructions are mutually inverse.
\end{void}

\begin{void}\label{geomgammanought}
Because over an algebraically closed field all the torsion points are
rational, such a base scheme is much easier to handle. Let
$\E/\overline{k}$ be an elliptic module over the algebraically closed
$\A$-field $\A\to\overline{k}$. If $\n$ is not divisible by the
characteristic $\q$ of $\overline{k}$, a $\Gamma_0(\n)$-structure is a
submodule $H$ of $\E[\n](\overline{k})$ that is isomorphic to $\A/\n$.
If $\n=\q^n$, and $\E$ is ordinary, then a
$\Gamma_0(\n)$-structure $H$ is either 
  $$ H = \E[\n]^\red $$
or the effective cartier divisor
  $$ H = \#(\A/\n)[0] $$
where $0$ is the zero section of $\GG_{a,\L}$. If $\E$ is super-singular,
the only $\Gamma_0(\n)$-structure is $ H = \#(\A/\n)[0] $.
\end{void}

\subsection{Mixed structures}

\begin{void}
Combinations of the above level structures will also be used. Let $\n$,
$\n_1$, and $\n_0$ be three ideals in $\A$. Giving a
$\Gamma(\n)\times\Gamma_1(\n_1)\times\Gamma_0(\n_0)$-structure on an
elliptic module $\E$ is the same as giving three level structures:
a $\Gamma(\n)$-structure, a $\Gamma_1(\n_1)$-structure, and a
$\Gamma_0(\n_0)$-structure.
\end{void}

\subsection{The action of $\GL(2, \A/\n)$}

\begin{void}
Let us fix an ideal $\n$ and an elliptic module $\E/S$.

The group $\GL(2,\A/\n)$ acts on the set of $\Gamma(\n)$-structures on $\E/S$
by left multiplication: an element
$m \in \GL(2,\A/\n)$ transforms a couple $(\E/S,\alpha)$ of a
module over $S$ and a $\Gamma(\n)$-structure on $\E$ as
 $$ (\E/S,\alpha) \mapsto (\E/S,\alpha \circ m)$$
and this action commutes with base change. Similarly, the group
$(\A/\n)^\times$ acts on $\Y_1(\n)$. Over an algebraically closed base field
$\overline{k}$, we can redefine a $\Gamma_0(\n)$-structure as a
$\Gamma_1(\n)$-structure modulo the left action of $(\A/\n)^\times$ on it.
\end{void}

\section{Modular schemes}

In this chapter we will represent functors that classify drinfeld modules
of rank 2 with certain level structures by affine schemes, and briefly mention
the canonical compactifications of them. We will only consider modules of rank
2, so-called elliptic modules, and as a consequence the fibers of the moduli
schemes over $\A$ will be curves.  

\subsection{The three basic moduli problems}

Compare to \cite{Katz85}. We do not consider balanced level structures
because we lack a good moduli interpretation in terms of dual isogenies, when
the ring $\A$ is not a principal ideal domain. The
associated moduli problem can however be defined as a quotient of the
full level moduli problem.

\begin{void}
Our three main moduli problems are functors on $(\Sch/\A)$, that associate
with a scheme the set of elliptic modules with certain level structure,
up to isomorphy. Concretely, we define, for an ideal $\n$ of $\A$ the
\emph{full level $\n$-structure moduli problem} or the 
\emph{$\gamma(\n)$-structure moduli problem}
  $$ \Y(\n): S \mapsto \Y(\n)(S) =
  \{\text{elliptic modules over } S \text{ with } \Gamma(\n)
  \text{-structure} \}/\cong $$
as well as the \emph{$\gamma_1(\n)$-structure moduli problem}
  $$ \Yi(\n): S \mapsto \Yi(\n)(S) = 
  \{\text{elliptic modules over} S \text{ with } \Gamma_1(\n)
  \text{-structure} \}/\cong $$
and the \emph{$\gamma_0(\n)$-structure moduli problem}
  $$ \Yo(\n): S \mapsto \Yo(\n)(S) = 
  \{\text{elliptic modules over } S \text{ with } \Gamma_0(\n)
  \text{-structure} \}/\cong $$
\end{void}

\begin{void}
Given an elliptic module $\E/S/\A$, we write, by abuse of notation, 
  $$ \Y(\n)(\E/S) =  \{\Gamma(\n)\text{-structures on } \E/S\} $$
This is in fact functorial by base change. In the same way we write
  $$ \Yi(\n)(\E/S) =  \{\Gamma_1(\n)\text{-structures on } \E/S\} $$
  $$ \Yo(\n)(\E/S) =  \{\Gamma_0(\n)\text{-structures on } \E/S\} $$
can define new functors $\Yi(\n)$ and 
$\Yo(\n)$. We call these functors on $(\text{elliptic modules})$
\emph{relative moduli functors}, to differentiate them from the
\emph{moduli functors} on $(\Sch/\A)$.
\end{void}

\begin{void}
For mixed structures we use the following notation. We simply concatenate
the symbols for functors. For example, the moduli problem
$\Gamma(\n)\times\Gamma_0(\m)$ defines a moduli functor denoted by
$\Y(\n)\Yo(\m)$.
\end{void}

\subsection{Relative representability}

In this section we will prove the following result, which we call
\emph{relative representability}.

\begin{theorem}\label{relrep}
  Given an elliptic module  $\E/S/\A$ and an ideal $\goth{n}$, then the
  three functors on $(\Sch/S)$
  $$ T \mapsto \Y(\n)(\E_T/T)=\{\Gamma(\n)\text{-structures on }\E_T\}$$
  $$ T \mapsto \Yi(\n)(\E_T/T)=\{\Gamma_1(\n)\text{-structures on }\E_T\}$$
  $$ T \mapsto \Yo(\n)(\E_T/T)=\{\Gamma_0(\n)\text{-structures on }\E_T\}$$
  are represented by finite affine schemes over $S$.
\end{theorem}

  The method is the same as used in \cite{Katz85}, in the
  context of elliptic curves for $\Gamma$- and $\Gamma_1$-structures. For
  $\Gamma_0$-structures we give a much easier proof, using faithfully flat
  descent.

\begin{lemma}
  For a fixed elliptic module $\E/S/\A$, the functors 
    $$ T \mapsto \Hom_\A(\A/\n,\E[\n](T))  $$
  and
    $$ T \mapsto \Hom_\A((\A/\n)^2,\E[\n](T)) $$
  are represented by the finite $S$-schemes $\E[\n]$ and
  $\E[\n] \times_S \E[\n]$ respectively.
\end{lemma}

\begin{proof}
   Because the $\A$-module $\E[\n](T)$ descends to an $\A/\n$-module
   (\ref{torsionscheme}), we have canonical isomorphisms:
     $$ \Hom_\A(\A/\n,\E[\n](T))
        = \Hom_{\A/\n}(\A/\n,\E[\n](T))
        = \E[\n](T) $$
   This proves the first claim, and with this, the second result follows
   from the canonical isomorphism
     $$ \Hom_\A((\A/\n)^2,\E[\n](T)) = 
        \left(\Hom_\A(\A/\n,\E[\n](T))\right)^2 $$
\end{proof}

\begin{lemma}\label{reli}
  The functor
    $$ T \mapsto \Yi(\n)(\E_T/T) $$
  is represented by a closed subscheme of $\E[\n]$.
\end{lemma}

\begin{proof}
  By definition, the considered functor is a subfunctor of the representable
  functor $T\mapsto\Hom_\A(\A/\n,\E[\n](T))$. An $\A$-morphism
   $$ \alpha: \A/\n \to \E[\n](T) $$
  is a $\Gamma_1(\n)$-structure on $\E_T$ if and only if it induces an
  inequality of effective cartier divisors
   $$ [\alpha(\A/\n)] \leq [\E_T[\n]]. $$
  This condition is represented by a closed subscheme of $\E[\n]$ by
  (\cite{Katz85} 3.1.4).
\end{proof}

\begin{lemma}\label{rel}
  The functor
  $$ T \mapsto \Y(\n)(\E_T/T) $$
  is represented by a closed subscheme of $\E[\n]\times_S\E[\n]$.
\end{lemma}

\begin{proof}
  The proof is essentially the same as the previous one. Given an
   $$ \alpha: (\A/\n)^2 \to \E[\n](T), $$
  it is a $\Gamma(\n)$-structure if and only if 
   $$ [\alpha((\A/\n)^2)] = [\E_T[\n]]. $$
  But since both sides automatically have the same degree, it suffices to
  require the inequality
   $$ [\alpha((\A/\n)^2)] \leq [\E_T[\n]]. $$
\end{proof}

The proof for $\Gamma_0(\n)$-structures is more subtle, we will use the
technique of faithfully flat descent (\cite{Grothendieck95}, \cite{SGA1} VIII).
It also motivates our definition.
A priori, more definitions seemed plausible, but the given 
definition gives rise to a relatively representable moduli problem, a
necessity to obtain good moduli schemes.

\begin{lemma}\label{relo}
  Let $Z/S$ be the scheme representing the functor
  $$ T \mapsto \Yi(\n)(\E_T/T), $$
  then the functor
  $$ T \mapsto \Yo(\n)(\E_T/T) $$
  is represented by the quotient scheme $Z/G$, where the group
  $G=(\A/\n)^\times$ is the group of units of the ring $\A/\n$, acting on
  $Z$ by its action on level structures.
\end{lemma}

\begin{proof}
Given a $\Gamma_0(\n)$-structure $H$ on $\E_T$, we do a finite faithfully flat
base-change $T'\to T$ such that $H_{T'}$ has a base point. This choice of
base point yields a $\Gamma_1(\n)$-structure on $\E_{T'}$, and hence a
$T'$-valued point of $Z$. By composition $Z\to Z/G$ we then
find a section $T' \to Z/G$ of the quotient. Now by faithfully
flat descent, this section descends to a $T$-valued point. Indeed, write
$T''=T'\times_T T'$ and consider the following commutative diagrams, indexed
by $i=1,2$:
\[\begin{CD}
     T''      @>i>>    Z_{T''}   @>>>   (Z/G)_{T''} \\
  @V{p_i}VV                @.               @VV{i}V   \\
     T'       @>>>     Z_{T'}    @>>>    (Z/G)_{T'} \\
   @VVV                  @.                @VVV     \\
     T         @.                 @.       (Z/G)    
\end{CD}\]
The two cases correspond to the two projections $T''=T'\times_T T' \to
T'$. For the middle
composite horizontal arrow to descend, we have to show that
the composite top level arrow does not depend on $i$ either. The induced
elements of $Z(T'')$ are $\Gamma_1(\n)$-structures on $\E_{T''}$ refining
the given $H_{T''}$ (the composite arrow $T''\to T$ is independent of $i$).
Hence they differ by an automorphism in $G=(\A/\n)^\times$, and project
to the same element of $(Z/G)(T'')$.

Conversely, let $T \to Z/(\A/\n)^\times$ be a morphism of $S$-schemes.
Again, a finite faithfully flat extension yields a $T' \to Z$, namely
$T'=T\times_{Z/G}Z$. This corresponds to a $\Gamma_1(\n)$-structure on
$\E_{T'}$, and this in turn gives us a $\Gamma_0(\n)$-structure $H'$ over
$T'$. The group $G$ works on $T'$ by its action on $Z$, and also on $H'/T'$ by
this action on $T'$. We get a composite $H' \to T' \to T$ which is invariant
under $G$ and hence factors over $H=H'/G \to T$. Hence we get a cartesian
square
\[\begin{CD}
    H'    @>>>     T' \\
  @VVV           @VVV \\
    H     @>>>     T
\end{CD}\]
and a canonically determined element of $\Yo(\n)(\E/T)$. 

The above constructions are each others inverses, and we have established a
canonical bijection, functorial in $T$, 
 $$ \Yo(\n)(\E_T/T) = \Hom(T,Z/G). $$

\end{proof}
 
\begin{theorem}\label{relquot}
 Fix an elliptic module $\E/S$. Let $Z(\n)$, $Z_1(\n)$ and $Z_0(\n)$ be the
 representing schemes for the relative moduli functors corresponding to
 $\Gamma(\n)$, $\Gamma_1(\n)$ and $\Gamma_0(\n)$. We have the following
 quotients by subgroups of $\GL(2,\A/\n)$, acting on level structures:
 $$ Z(\A)=Z(\n)/\left(
   \begin{array}{cc} \ast & \ast \\
                     \ast & \ast \\
   \end{array}\right)$$
 $$ Z_0(\n)=Z(\n)/\left(
   \begin{array}{cc} \ast & \ast \\
                       0  & \ast \\
   \end{array}\right)$$
 $$ Z_1(\n)=Z(\n)/\left(
   \begin{array}{cc}   1 & \ast \\
                       0 & \ast \\
   \end{array}\right)$$
 If $\m \subset \n$ is an inclusion of ideals, we also have:
 $$ Z(\n)=Z(\m)/G $$
 where $G$ is the kernel of the reduction map $\GL(2,\A/\m) \to
\GL(2,\A/\n)$.
\end{theorem}

\begin{proof}
  All these claims can be proven in precisely the same way as the second one,
  about $Z_0$ (\ref{relo}).
\end{proof}

\subsection{Coarse and fine moduli schemes}

A pioneering book on the topic of moduli spaces is ``Geometric Invariant
Theory'' (\cite{Mumford94}).

\begin{void}
Let $\cal{F}$ be a functor from $(\Sch/S)$ to $(\text{sets})$. A
\emph{coarse moduli scheme} for $\cal{F}$ is a couple $(Y,\Phi)$ of an
$S$-scheme $Y$ and a morphism of functors
  $$ \Phi: \Hom(Y,-) \to \cal{F}(-) $$
such that
\begin{enumerate}
  \item $\Phi(\overline{k}): Y(\overline{k}) \to \cal{F}(\overline{k})$
     is a bijection for every algebraically closed field $\overline{k}$
  \item $(Y,\Phi)$ satisfies the following universal property: if $Y'$ is a
     scheme over $S$, and $\Phi': \Hom(Y',-) \to \cal{F}(-)$ a morphism of
     functors, there is a unique morphism of functors
     $f: \Hom(Y,-) \to \Hom(X,-)$ such that $\Phi' = f \circ \Phi$.
\end{enumerate}

If, moreover, $\Phi$ is a functorial isomorphism, we call $(Y,\Phi)$ a
\emph{fine moduli scheme}. In this case, the scheme $Y$ represents the
functor $\cal{F}$. Fine moduli schemes are a fortiori coarse moduli
schemes. From the universal property in the definition, it follows formally
that coarse moduli schemes are unique up to a unique isomorphism. 

Clearly a representable functor has a fine moduli scheme. A
non-representable functor does not always have a coarse moduli scheme.
\end{void}

\begin{definition}
  An ideal $\n$ of $\A$ is called \emph{admissible} if its factorization into
  prime ideals contains at least two different prime factors.
\end{definition}

\begin{theorem}\label{representability}
  Let $\n$ be an admissible ideal in $\A$ with at least two different prime
  factors. The moduli problem $\Y(\n)$ is representable by an affine
  $\A$-scheme. We denote it by $Y(\n)$.
  Moreover, for all ideals $\m$, $\Y(\n)\Yo(\m)$ and $\Y(\n)\Yi(\m)$ are 
  represented by affine $\A$-schemes, which
  we denote by $Y(\n)Y_0(\m)$ and $Y(\n)Y_1(\m)$.
\end{theorem}

This means that the scheme $Y(\n)$ is equipped with a \emph{universal
elliptic module} and a $\Gamma(\n)$-structure on it, such that every
elliptic module with such a level structure over a scheme $S$ is obtained
from this universal module by a base change $S \to Y(\n)$.
    
\begin{proof}
  (\cite{Drinfeld74E}, \cite{Saidi96} 3.2.6) 
  Assume $\p$ and $\q$ are different prime divisors of $\n$. We will first
  prove representability of $\Y(\n)$ over $\A(\p)$ and of $\A(\q)$, and then
  glue over $\A(\p\q)$. The result for the two mixed moduli problems
  follows by applying relative representability (\ref{reli} and \ref{relo})
  to the universal elliptic module $\E/Y(\n)$.

  Now let us proceed to prove that $\Y(\p^n)$, restricted to schemes over
  $\A(\p)$ is representable. The ring $\A$ can be presented in terms of
  generators and relations as
  $$ \A = \FF_q[a_1, a_2, \ldots a_m]/(f_1, f_2, \ldots ,f_{m-1}) $$
  We are going to construct an $\A$-algebra $B$, equipped with an
  elliptic module $\E=(\GG_a,\phi)$ for the trivial line bundle on $B$ and a
  $\Gamma(\n)$-structure on it. Then we will show that it is universal for
  the moduli problem restricted to $\A(\p)$-schemes. We construct $B$ in
  several steps, the construction is straightforward but involves a lot of
  generators and relations.
  \begin{enumerate}
   \item We add the coefficients $c_{i,j}$ of a generic elliptic module
   $\phi$ (in standard form) to $\A$:
    $$ \phi(a_i)=\sum{c_{i,j}\tau^j} $$
   \item We add the relations between the $c_{i,j}$ expressing that $\phi$
    is a ring homomorphism. Moreover we invert the coefficients of highest
    degree $c_{i,2\deg(a_i)}$, for all $i$.
   \item A $\Gamma(\p^n)$-structure on $\E$ is defined by 
    $$ \alpha : (\A/\p^n)^2 \to B : x \mapsto \alpha(x), $$
    hence we add the images $\alpha(x)$ to our set of generators.
   \item They must satisfy algebraic relations expressing that they are
    additive in $x$, and that they are $\A$-linear:
     $$ \phi_{a_i}\alpha(x)=\alpha(a_i x) $$
     $$ \alpha(\lambda x) = \lambda \alpha(x), \forall \lambda \in \FF_q $$
    Moreover, the image of $\alpha$ satisfies
     $$ \prod_{x \in (\A/\p^n)^2} (z-\alpha(x)) =  P(z) $$
    where $P$ is the monic polynomial corresponding to the $\q$-torsion of
    $\E$. We finish by adding all these relations.
  \end{enumerate}
  
  It follows from the construction that every elliptic module over an
  $\A$-algebra in which $\p$ is invertible can be obtained from $\E$ by a
  base change. Moreover, this base change is unique since our construction
  of $\E$ is rigid. Remark that if we allow algebras with
  characteristic intersecting $\p$ then we lose this uniqueness. From this
  universal property it follows immediately that $B$ represents the moduli
  problem over $\A(\p)$-algebras and even over $\A(\p)$-schemes. Applying
  relative representability to the universal elliptic module $\E/B$ shows
  that also $\A(\n)$, with $\p$ dividing $\n$, is representable over
  $\A(\p)$. 

  So assume now that $\n$ has two different prime divisors $\p$ and $\q$.
  We find that $\Y(\n)$ restricted to $\A(\p)$- or $\A(\q)$-schemes is
  represented  by some affine schemes $Y_1$ and $Y_2$ respectively. Since
  they represent the same functor over $\A(\p\q)$-schemes, there is a
  canonical isomorphism between open subschemes
    $Y_1 \otimes_\A \A(\p) \isomto Y_2 \otimes_\A \A(\q) $
  We find that $\Y(\n)$ is represented by the scheme $Y(\n)$ obtained by
  glueing $Y_1$ and $Y_2$ by this isomorphism.
\end{proof}

We now want to prove the existence of coarse moduli schemes for the other
moduli problems we defined. 

\begin{void}
Let $\n$ be any ideal in $\A$, and fix an admissible ideal $\m$ coprime to
$\n$.    
We define quotient schemes of fine moduli schemes by the finite group
$G=\GL(2,\A/\m)$, acting on $\Gamma(\m)$-structures.

 $$ Y(\n)=Y(\m\n)/G $$

 $$ Y_1(\n)=Y(\m)Y_1(\n)/G $$

 $$ Y_0(\n)=Y(\m)Y_0(\n)/G $$

These quotient schemes are
independent of the choice of $\m$, since the different quotient morphisms
coincide on geometric points. For the same reason, these definitions
are consistent with the ones already given for $Y(\n)$ with admissible
$\n$.
\end{void}

\begin{theorem}
  $Y(\n)$, $Y_1(\n)$ and $Y_0(\n)$ are coarse moduli schemes for the moduli
  problems $\Y(\n)$, $\Yi(\n)$ and $\Yo(\n)$, respectively.
\end{theorem}

\begin{proof}
  Let $Y$ be one of the schemes $Y(\n)$, $Y_1(\n)$ or $Y_0(\n)$, and let
  $\Y$ be the corresponding moduli functor ($\Y(\n)$, $\Yi(\n)$ or
  $\Yo(\n)$). By definition
    $$ Y = Y'/G $$
  for some finite group $G=\GL(2,\A/\m)$ and fine moduli scheme $Y'$,
  representing the moduli problem obtained by adding some
  $\Gamma(\m)$-structure to the old one.
  We will silently use that this is independent of the choice of $\m$.
  We proceed in several steps.

  {\sc Step 1.} We construct a morphism of functors
    $$ \Phi : \Y \to \Hom(-,Y) $$
  Fix an $\A$-scheme $S$. An element of $\Y(S)$ corresponds to an elliptic
  module $\E$ equipped with a certain level structure. 
  After a finite faithfully flat base change $T\to S$, we can chose a
  $\Gamma(\m)$-structure on $\E_T$. This induces a point $T\to Y_T$, which
  by faithfully flat descent descends to a point $S \to Y$.
  
  {\sc Step 2.} Remark that $ \Phi(\o{k}) $ is a bijection, for all
  algebraically closed $\o{k}$. Indeed, $Y(\o{k})=Y'(\o{k})/G$, and this is
  precisely $\Y(\o{k})$.

  {\sc Step 3.} (compare \cite{Mumford94} 5.2.4) There is a canonical
  bijection
    $$ \{\text{morphisms of functors } \phi:\Y \to \Hom(-,S) \}
       =
       \{\text{morphisms } f: Y \to S \}, $$
  functorial in $S$. We can see this as follows. The fine moduli scheme
  $Y'$ carries a universal elliptic module $\E^{\text{univ}}/Y'$ with some
  level structure. Forgetting the $\Gamma(\m)$-structure, we find an
  element $P \in \Y(Y')$. Now assume given a $\phi$ as above. Applying it
  to $P$ yields a morphism $Y' \to S$. But by construction it is invariant
  under the action of $\GL(2,\A/\m)$ and hence factors over a unique
  $f:Y \to S$. On the other hand, given an $f: Y \to S$, we can compose
  with $ \Phi : \Y \to \Hom(-,Y) $ to get a $\phi:\Y \to \Hom(-,S)$:
    $$ \Y \to \Hom(-,Y) \to \Hom(-,S) $$
  The correspondences $\phi \mapsto f$ and $f \mapsto \phi$ are each others
  inverses and we established the claimed functorial bijection.

  {\sc Step 4.} Combining the previous steps we can verify that $Y$ is
  a coarse moduli scheme.
\end{proof}

\begin{void}
To study the coarse moduli schemes, one often covers them first with a fine
moduli scheme by adding, for example, some $\Gamma(\p\q)$-structure, for
suitable primes $\p$ and $\q$ to the moduli problem. The new, possibly
mixed, moduli functor is representable. We can study the properties of
the new, fine, moduli scheme by studying the functor it represents. Then we
descend to the original, coarse, moduli scheme by taking the quotient by a
suitable group. This procedure can often be avoided by using the
machinery of \emph{stacks} as described in \cite{Deligne69}.
\end{void}

\subsection{Deformation theory}

In this section we will study the regularity and dimensionality of the
fine and coarse moduli schemes, and of their fibers. To do so, we first
study regularity in the case of fine moduli schemes. And then take
quotients by suitable finite groups, acting without fix-points, to obtain
results for the coarse moduli schemes. As everywhere else in this chapter,
we restrict to the case of rank two modules, but the methods used in this
section are easily generalized to higher rank drinfeld modules.

\begin{void}
In what follows, $k$ is an algebraically closed field. $k[\epsilon]$ is the
$k$-algebra determined by the relation $\epsilon^2=0$. An algebraic
structure $\o{A}$ over $k[\epsilon]$ is a \emph{lift} of the structure $A$ over
$k$, if $A$ is obtained from $\o{A}$ by base change via
 $$ k[\epsilon]\to k: \epsilon \mapsto 0 $$
The lift is said to be \emph{trivial} if it is obtained from $A$ by base
change via the inclusion $k \injto k[\epsilon]$.
\end{void}

\begin{lemma}
  Assume $k$ has the structure of an $\A$-algebra by a non-zero morphism
  $\gamma:\A\to k$. Then the space of lifts of $\gamma$ to 
   $$ \o{\gamma}:\A\to k[\epsilon] $$
  is one-dimensional over $k$.
\end{lemma}

\begin{proof}
  This is just a restatement of the regularity of the curve $\A$ over
  $\FF_q$.
\end{proof}

\begin{proposition}
  Let $\E$ be an elliptic module over $\gamma:\A\to k$. Given a fixed lift
  $\o{\gamma}:\A \to k[\epsilon]$ of $\gamma$, the $k$-space of
  lifts $(\o{\gamma},\o{\E}/k[\epsilon])$ of $(\gamma,\E/k)$ is
  one-dimensional.
\end{proposition}

\begin{proof}
  See also (\cite{Laumon96} 1.5) and (\cite{Dokchitser00} 2.2).
  The following proof is based on the proof in \cite{Laumon96}.
  It was simplified by Marius van der Put.
  
  Remark that all line bundles on $k[\epsilon]$ are trivial.

  Let us first lift $\E$ without fixing a given lift of $\gamma$.

  The elliptic module $\E$ is given by a morphism $\phi:\A\to k\{\tau\}$ of
  rings. Now consider a lift
   $$\o{\phi}=\phi+\epsilon D :\A\to k[\epsilon]\{\tau\}$$
  of $\phi$. Demanding that $\o{\phi}$ is a ring homomorphism, implies that
  $D$ is a derivation $D:A\to M$ to the bimodule $M=k\epsilon\{\tau\}$, on
  which $a\otimes b \in \A\otimes\A$ acts as
   $$ amb = \phi(a)m\phi(b) = \gamma(a)m\phi(b). $$
  Conversely, every such derivation yields a ring homomorphism
  $\o{\phi}=\phi+\epsilon D$. 

  A lift $\o{\phi}=\phi+\epsilon D$ is isomorphic to the trivial lift
  precisely when there is a fixed $m\in M$ such that for all $a$,
    $$ D:A\to M: a \mapsto am-ma. $$

  Now we introduce the \emph{universal bimodule derivation} (\cite{EGA0} \S20).
  We write $I$ for the kernel of the multiplication map
    $$ m: \A\otimes_{\FF_q}\A \to \A : a\otimes b \mapsto ab $$
  and equip $I$ with the obvious $\A/\FF_q$-bimodule structure. The bimodule
  derivation
    $$ d: \A \to I : a \mapsto a\otimes1-1\otimes a $$
  has the following universal property. Every derivation $D:\A\to N$ of $\A$
  in an $\A/\FF_q$-bimodule $N$ factors in a unique way as $D = l \circ d$,
  where $l:I\to N$ is a bimodule homomorphism. 
  
  In our situation, we find the space of lifts of ring homomorphisms is
  $\Hom_{\A\otimes\A}(I,M)$, and the subspace of lifts equivalent to the
  trivial one $\Hom_{\A\otimes\A}(\A\otimes\A,M)$. Hence the sought
  deformation space $\text{Def}$ is determined by the exact sequence
    $$ \Hom_{\A\otimes\A}(\A\otimes\A,M) \to
       \Hom_{\A\otimes\A}(I,M) \to \text{Def} \to 0 $$

  Let $l:I\to\A\otimes\A$ be an $\A\otimes\A$-morphism. The composition with
  the multiplication map $m$ is a morphism $I\to \A$. Since this morphism is
  zero on $I^2$, it factors over a morphism $s(l): I/I^2 \to \A$. We have
  constructed a map
  $s:\Hom_{\A\otimes\A}(I,\A\otimes\A)\to\Hom_{\A\otimes\A}(I/I^2,\A)$.
  To calculate its kernel, note that $s(l)=0$ is equivalent to $l(I)\in I$.
  The ideal $I\subset A\otimes A$ corresponds geometrically
  to the diagonal divisor on the surface $\Spec(A\otimes A$), hence
  $I$ is is a projective $\A\otimes\A$ module of rank $1$,
  and thus the isomorphism
  $\Hom_{\A\otimes\A}(I,I)\isomto\A\otimes\A$ holds.
  Therefore, $l(I)\in I$
  is equivalent to the existence of a continuation of $l:I\to\A\otimes\A$ to
  an $\A\otimes\A$-morphism $\A\otimes\A\to\A\otimes\A$. We have proven the
  exactness of the sequence of $\A\otimes\A$-modules
   $$ \Hom(\A\otimes\A,\A\otimes\A)\to\Hom(I,\A\otimes\A)
      \to\Hom(I/I^2,\A)\to 0 $$
  in which the middle map is the $s$ constructed above. The
  $\A\otimes\A$-structure on $\A$ is induced by the multiplication morphism
  $m$. Tensoring this exact sequence with $M$ yields the previous exact
  sequence and we find
   $$ \text{Def} \cong \Hom_{\A\otimes\A}(I/I^2,\A)\otimes_{\A\otimes\A} M $$
 
  The ideal $I=\ker( m:\A\otimes\A\to\A)$ acts as $0$ on both $I/I^2$ and
  $\A$, hence $\Hom_{\A\otimes\A}(I/I^2,\A)=\Hom_\A(I/I^2,\A)$ 
  which is projective of rank $1$ over $A$ since it is dual 
  to the module of relative differentials of $\A$ over $\FF_q$.
  
  Now the following lemma allows us to conclude that
  $\text{Def}\cong k^2$ as $k$-vector space. 
  
  Hence, we have proven that the space of lifts of $\E$ to arbitrary lifts
  of $\gamma$ is two-dimensional. We still have to check that it projects
  surjectively to the one-dimensional space of lifts of $\gamma$. This
  follows from the above proof, by tracking in every step what happens to
  the lifts of $\gamma$ induced by the lifts of $\E$.
\end{proof}

\begin{lemma}
  Let $\gamma:\A\to k$ be an $\A$-field, and $\E=(\GG_{a,k},\phi)$ be an
  elliptic module over it. The $k\otimes_{\FF_q}\A$-module $k\{\tau\}$ given
  by the action
    $$ (\lambda\otimes a)\alpha = \lambda \alpha \phi_a \in k\{\tau\} $$
  is locally free of constant rank 2.
\end{lemma}

\begin{proof}
  First we change the base field of the curve $\A$ by $\FF_q \injto k$.
  We write $\A_k=\A\otimes_{\FF_q}k$ and $\gamma_k=\gamma\otimes_{\FF_q}k$.
  Choose a non-constant element $a\in\A$ of degree $d$. Then $\A_k$ is
  locally free of rank $d$ over the polynomial ring $k[a]\injto\A_k$.
  Moreover, by using the left euclidean algorithm we can write every
  element of $k\{\tau\}$ uniquely as a $k$-linear combination of terms
    $$ \tau^n\phi_{a^m}=\tau^n\phi_a^m $$
  where $0\leq n < 2d$. Hence $k\{\tau\}$ is free of rank $2d$ over $k[a]$.
  Because $\FF_q$ is precisely the field of constants of $\A$, $\A_k$ is
  integral and the result follows.
\end{proof}

\begin{proposition}
  Assume $\n$ is an ideal of $\A$ which is coprime to the characteristic
  $\gamma:\A\to k$. Let $(\o{\gamma},\o{\E})$ be a lift of $(\gamma,\E)$.
  A $\Gamma(\n)$, $\Gamma_1(\n)$ or $\Gamma_0(\n)$-structure lifts uniquely
  to a similar structure on $\o{\E}$.
\end{proposition}

\begin{proof}
  The $\n$-torsion of $\E$ is given by a separable, additive polynomial $f$
  over $k$. This implies the coefficient of its linear term is non-zero,
  say $c$. Similarly, $\o{\E}[\n]$ is given by a lift $\o{f}$ of $f$ to
  $k[\epsilon]$. Now let $u\in k$ be a zero of $f$, which means it
  corresponds to a $k$-rational torsion point of $\E$. For the lift
  $u + \epsilon v \in k[\epsilon]$ to be a $k[\epsilon]$-rational torsion
  point of $\o{\E}$, it suffices that $\o{f}(u+\epsilon v)=0$. We
  calculate 
    $$ \o{f}(u+\epsilon v)= \o{f}(u) + \epsilon vc $$
  where $\o{f}(u)$ is an element of the radical $k\epsilon$ of
  $k[\epsilon]$. Hence, since $c \neq 0$, we find a unique $v\in k$ such
  that $\o{f}(u+\epsilon v)=0$. This already proves the proposition for
  $\Gamma(\n)$ and $\Gamma_1(\n)$-structures. 

  Now for $\Gamma_0(\n)$-structures, remark that we have just shown that
  both $\E[\n](k)$ and $\o{\E}[\n](k[\epsilon])$ are isomorphic to
  $(\A/\n)^2$. Therefore choosing a $\Gamma_0(\n)$-structure on one of them
  is nothing but choosing a submodule $\A/\n \subset (\A/\n)^2$ and no base
  extensions are needed to make torsion rational. Now the lemma follows
  immediately by the lifting bijection from $\E[\n](k)$ to
  $\o{\E}[\n](k[\epsilon])$ we established.
\end{proof}

\begin{corollary}\label{dim}
  For all $\n$, the moduli schemes $Y(\n)$, $Y_1(\n)$ and $Y_0(\n)$ have
  dimension $2$ over $\FF_q$, the field of constants of $\A$.
\end{corollary}

\begin{proof}
  Let $\p$ be a prime ideal. It follows from the proof of
  \ref{representability} that $\Y(\p)$ is representable by an affine scheme
  $Y$ over $\A(\p)$. Moreover, it follows from the given construction that this
  affine scheme is reduced. The previous proposition calculates the tangent
  space at a geometric point of $Y$, and since it is everywhere
  2-dimensional, $Y$ must be
  regular of dimension $2$ over $\FF_q$. Moreoever, $Y$ must coincide with
  the restriction $Y(\p)\otimes_{\A} \A(\p)$ of the coarse moduli scheme
  $Y(\p)$, which is therefore also 2-dimensional. All the other moduli
  schemes are finite covers of $Y(\p)$, or finite quotients of these finite
  covers, and hence they are 2-dimensional.
\end{proof}

\begin{proposition}
  Let $\p$ be the kernel of $\gamma$, and assume it is non-zero. Let $\E$
  be an elliptic module over $\gamma:\A\to k$, equipped with a
  $\Gamma(\p^s)$-structure $\alpha$. The lifted couple $(\o{\gamma},\o{\E})$
  allows a lift $\o{\alpha}$ of $\alpha$ if and only if it is the trivial
  lift of $(\gamma,\E)$. In this case, the space of lifts of $\alpha$ is
  two-dimensional over $k$.
\end{proposition}

\begin{proof}
  Write $\n=\p^s$. The $\Gamma(\n)$-structure
    $$ \alpha: (\A/\n)^2 \to \GG_{a,k}(k)=k $$
  is given by its values on a free basis of $(\A/\n)^2$ as $\A/\n$-module.
  We set $u_1=\alpha(1,0)$, $u_2=\alpha(0,1)$, and similarly
  $\o{u_1}=u_1+\epsilon v_1=\o{\alpha}(1,0)$ and
  $\o{u_2}=u_2+\epsilon v_2=\o{\alpha}(0,1)$ for the lift
    $$ \o{\alpha}:(\A/\n)^2\to\GG_{a,k[\epsilon]}(k[\epsilon])=k[\epsilon] $$
  of $\alpha$.
  
  We must have that the $\n$-torsion schemes
  (divisors) of $\E$ and $\o{\E}$ correspond to the polynomials
    $$ f = \prod_{x\in(\A/\n)^2}(X-\alpha(x)) \text{ and }
      \o{f} = \prod_{x\in(\A/\n)^2}(X-\o{\alpha}(x)) $$
  Because $\n$ is a power of the characteristic, the kernel of $\alpha$ is
  non-trivial and therefore contains a submodule acting as $\FF_q$. Now
  using the formula
    $$ \prod_{\lambda\in\FF_q}(X-u-\lambda v \epsilon) = (X-u)^q $$
  we find that $\o{f}=f$. 

  Now chose $t>0$ such that $\n^t$ is a principal ideal, say
  generated by $a$. Then the $\n^t$-torsion of $\E$ is defined by the
  polynomial $f(f(\cdots(f)))$ ($t$ iterations), and the $\n^t$-torsion of
  $\o{\E}$ by $\o{f}(\o{f}(\cdots(\o{f})))$. Since these are equal, we find
    $$ \phi_a = \o{\phi}_a $$
  Because $\o{\phi}_a$ has positive height, we also derive from
    $$ \o{\phi}_b \o{\phi}_a = \o{\phi}_a \o{\phi}_b $$
  that $\o{\phi}_b = \phi_b$ for every $b\in\A$. Hence, $\o{\gamma}$ and
  $\o{\E}$ must be trivial lifts. Now on the trivial lift, we have a
  two-dimensional space of candidate lifts $\o{\alpha}$, namely every
  candidate is determined by choosing $v_1$ and $v_2$. Every such choice
  yields a $\Gamma(\n)$-structure on $\o{\E}$ because of the same formula
  that was used earlier in this proof.
\end{proof}

\begin{proposition}
  Let $\p$ be the kernel of $\gamma:\A\to k$ and assume given an elliptic
  module over $\gamma:\A\to k$ with a $\Gamma_1(\p^s)$-structure $\alpha$.
  We want to understand the structure of the space of lifts
  $(\o{\gamma},\o{\E},\o{H})$ of $(\gamma,\E,H)$.
  We consider two cases:
  \begin{enumerate}
    \item If the image of $\alpha$ is \'{e}tale ($\E$ ordinary), every choice
          of $\o{\gamma}$ allows a unique choice of $\o{\E}$ which, in turn,
          has a one-dimensional space of compatible lifts of $\alpha$.
    \item If $\alpha$ is local ($\E$ ordinary or supersingular), every choice
          of $(\o{\gamma},\o{\E})$ admits a one-dimensional space of
          compatible lifts of $\alpha$.
  \end{enumerate}
\end{proposition}

\begin{proof}
Write $\n=\p^s$. The $\Gamma_1(\n)$-structure
 $$ \alpha : (\A/\n) \to \GG_{a,k}(k)=k $$
is determined by the image of $1 \in \A/\n$. We denote it by $u\in k$.
The lift $\o{\alpha}$ of $\alpha$ is determined by $\o{u} = u + \epsilon
v$. Let $f$ be the monic polynomial corresponding to the $\n$-torsion divisor,
that is, let $f$ be such that
 $$ \E[\n] = k[X]/(f(X)) $$
If we write $d$ for $\deg(\n)$, we know that $f$ has the form
 $$ f(X) = a_d X^{q^d} + a_{d+1} X^{q^{d+1}} + \ldots + X^{q^{2d}} $$
and that $\alpha$ is a $\Gamma_1(\n)$-structure is expressed by $f(u)=0$.
The $\n$-torsion of a lift $\o{\E}$ corresponds to a lift of $f$:
 $$ \o{f} = f + \epsilon g $$
Now $\o{u}=u+\epsilon v$ corresponds to a $\Gamma_1(\n)$-structure if
$\o{f}(\o{u})=0$. Since $f'=0$, this translates to $g(u)=0$. 

We conclude that if $u=0$, that is $\alpha$ is a purely local level
structure, every lift of $(\gamma,\E)$ allows a one-dimensional space of
lifts of $\alpha$.

On the other hand, if $u\neq0$, we find that for every lift of $\gamma$, a
unique compatible lift of $\E$ allows for lifts of $\alpha$. Moreover, this
space of lifts of $\alpha$ is one-dimensional.
\end{proof}

\begin{proposition}
  Let $\p$ be the kernel of $\gamma:\A\to k$ and assume given an elliptic
  module over $\gamma:\A\to k$ with a $\Gamma_0(\p)$-structure $H$.
  We want to understand the structure of the space of lifts
  $(\o{\gamma},\o{\E},\o{H})$ of $(\gamma,\E,H)$.
  As we know, there are two possibilities for $H$ when $\E$ is ordinary,
  and one when $\E$ is super-singular.
  \begin{enumerate}
    \item If $H=\E[\p]^\red$ ($\E$ ordinary), every choice of
          $\o{\gamma}$ allows a unique choice of $\o{\E}$ which, in turn,
          has a one-dimensional space of compatible lifts of $H$.
    \item If $H=\E[\p]^0$ ($\E$ ordinary), every choice of
          $(\o{\gamma},\o{\E})$ admits a unique compatible lift of $H$
    \item If $H$ is the unique $\Gamma_0(\p)$-structure on the
          super-singular elliptic module $\E$, only the trivial lift
          $\o{\gamma}$ of $\gamma$ has elliptic modules that allow lifting
          of the level structure. Moreover, every lift of $\E$ to an
          $\o{\E}$ over $\o{\gamma}:\A\to k[\epsilon]$ allows a
          one-dimensional space of lifts of $H$.
  \end{enumerate}
\end{proposition}

\begin{proof}
 Let $f(X)$ be the monic $\p$-torsion polynomial, and $f(X) + \epsilon g(X)$
 be a lift of it. A given $\Gamma_0(\p)$ structure corresponds to a monic
 polynomial $h(X)$, dividing $f(X)$. The proposition is now a purely
 algebraic exercise: express that a \emph{monic} lift of $h$ divides the
 given lift of $g$. The reasoning is very similar to that of the previous
 proposition. 
\end{proof}  

\begin{remark}
 It appears to be more difficult to examine higher $\p$-powers, although the
 same techniques should work.
\end{remark}

\begin{theorem}\label{smoothness}
  The surfaces $Y(\n)$ and $Y_1(\n)$, over the field of constants $\FF_q$
  of the curve $\A$ are regular. The surface $Y_0(\n)$ is regular outside
  super-singular points in fibers above $\p$'s for which $\p^2$ is a factor
  of $\n$.
  
  The curves $Y(\n)\otimes \kk(\p)$, $Y_1(\n)\otimes \kk(\p)$ and
  $Y_0(\n)\otimes \kk(\p)$ are regular if $\n$ and $\p$ are coprime.
  When $\p$ is a factor of $\n$, $Y(\n)\otimes \kk(\p)$ is not reduced, but
  $Y_1(\n)\otimes \kk(\p)$ is reduced and regular. 
  
  When $\n$ and $\p$ are coprime, the curve $Y_0(\n\p)\otimes\kk(\p)$ is
  regular outside the supersingular points. 
\end{theorem}

\subsection{Compactification}

The fine and coarse moduli schemes introduced in the previous section fail
to be proper. They can, however, be embedded as open dense subschemes of
proper schemes over $\A$. In the case of representable moduli problems, and
thus fine moduli schemes, the result is the following, due to Drinfel'd.

\begin{theorem}\label{compactification}
  Let $\n$ be an admissible ideal. There exists a
  unique regular scheme $X(\n)$ satisfying:
  \begin{enumerate}
    \item $Y(\n)$ is an open dense subscheme of $X(\n)$,
    \item $X(\n)$ is proper over $\A$,
    \item $X(\n)-Y(\n)$ is finite over $\A$,
    \item $X(\n)$ is smooth of relative dimension 1 over $\A(\n)$,
    \item $X(\n)$ is smooth of relative dimension 2 over $\FF_q$.
  \end{enumerate}
\end{theorem}

\begin{proof} See \cite{Lehmkuhl00}. The theorem already appears in Drinfel'd's paper
\cite{Drinfeld74E}.
\end{proof}

By taking quotients we can extend this theorem to the coarse moduli schemes.
We call the resulting $X(\n)$, $X_1(\n)$ and $X_0(\n)$ the
\emph{compactifications} of respectively $Y(\n)$, $Y_1(\n)$, $Y_0(\n)$.
We have the following result about regularity at cusps:

\begin{theorem}\label{allsmooth}
 Let $n$ be any ideal, then $X(\n)$, $X_1(\n)$ and $X_0(\n)$ are smooth
 curves over $\A(\n)$ and smooth surfaces over $\FF_q$.  
\end{theorem}
  
\begin{proof}
 On the affine parts $Y_\ast(\n)\injto X_\ast(\n)$ the claimed smoothness
 has been established already (see \ref{smoothness}). It remains
 to extend it to the compactification. Let $\m$ be a principal
 ideal contained in $\n$, and assume that $\m$ has at least
 two different prime divisors.
 In \cite{vanderHeiden03} we find an explicit description of the
 formal neighbourhood of cusps
 of $X(\m)$ above $\A(\m)$. It is the formal spectrum of a ring of
 the form
  $$ C = \oplus R[[X]] $$
 where $R$ is an integrally closed domain, finite over $\A(\m)$,
 and $\oplus$ a finite direct sum.
 The action of $\GL(2,\A/\m)$ is simultaneously by permutations of
 the components, and by normalized $R$-linear
 automorphisms of the $R[[X]]$
 (those that send $X$ to a power series of the form
 $X+\text{higher order terms}$). Now consider a subgroup $G$ of
 $\GL(2,\A/\m)$, acting stabily on one of the components,
 then it follows from the following lemma that the ring of
 invariants is also of the general form $\oplus R[[X]]$.
 Since the curves of the proposition can be obtained
 as quotients by such $G$ of $X(\m)$, it follows that
 they are smooth at the cusps.
\end{proof}

\begin{lemma}
 Let $R$ be a Dedekind ring and $G$ a finite group
 acting faithfully on $R[[X]]$, such that for ever $g\in G$
 $g(X)=uX$ with $u$ a unit in $R[[X]]$. Then
 $R[[X]]^G=R[[Y]]$, with $Y=\prod_{g\in G} g(Y)$.
\end{lemma}

\begin{proof}
 Clearly $R[[X]]^G\supset R[[Y]]$. Now let $s$ be a place
 of $R$ and denote by $R_s$ the local ring at $s$ and by
 $Q_s$ it's fraction field. Then the fraction field
 of $R_s[[X]]$ is $Q_s((X))$. We have that
 $Q_s((X))^G=Q_s((Y))$. Now since both $R_s[[X]]^G$ and
 $R_s[[Y]]$ are integrally closed in this field,
 and since one contains the other we have that
 $R_s[[X]]^G=R_s[[Y]]$ for every $s$, and hence that
 $R[[X]]^G=R[[Y]]$.
\end{proof}

\section{Curves of increasing genus having many points}

Fix a prime power $q$. It has been shown by Drinfel'd and
Vl\v{a}du\c{t} that a curve of genus $g$ over $\FF_q$ can have
at most
 $$ (\sqrt{q} - 1 + o(1))g $$
rational points, when $g\to \infty$ (see \cite{Vladut83}).
An infinite sequence of curves $C_i$ of increasing genus is said to be
\emph{asymptotically optimal} if $\# C_i(\FF_q)/g(C_i)$ tends
to $\sqrt{q}-1$.

The above bound is sharp at least when $q$ is a square, since several
families of (reductions of) modular curves, shimura curves and
drinfeld modular 
curves are known to be assymptotically optimal over a quadratic
extension of their finite field of definition. 
A paper by Noam Elkies (\cite{Elkies01}) gives explicit equations of a few towers of reductions of drinfeld modular curves, and shows that they are assymptotically optimal over the quadratic extension
of their field of definition. Moreover, the paper suggests
that assymptotic
optimality can be proven for more general drinfeld moduler curves in the same way Yasutaka Ihara proved a similar result for classical modular curves. Using the necessary tools from
the previous chapters, we will formulate and prove such a result here.
We note that the result itself already appears in \cite{Vladut91},
but a correct proof over general rings $\A$ seems to be missing
from the literature.

Fix a principal non-zero prime ideal $\p$ in $\A$ and an ideal
$\n\subset \A$ coprime to $\p$. Denote the degree of $\p$ by $m$, so
$\kk(\p)\cong \FF_{q^m}$.  The strategy, due to Ihara, is as
follows. Using
the modular interpretations, we define a curve $T$ on the surface
$X_0(\n)\times X_0(\n)$ so that the reduction of $T$ modulo $\p$
is supported at the union of the graphs of the
$q^m$-th power frobenius morphisms
$X_0(\n)\otimes\kk(\p)\to X_0(\n)\otimes\kk(\p)$ and
$X_0(\n)\otimes\kk(\p)\leftarrow X_0(\n)\otimes\kk(\p)$. These
two graphs intersect above the $\FF_{q^{2m}}$-rational points of $X_0(\n)\otimes\kk(\p)$. Hence the number of
singular points of the reduction of $T$ modulo $\p$
gives a lower bound for the number of rational points of
$X_0(\n)\otimes\kk(\p)$. Also, the hurwitz formula gives an
upper bound on the genus of $X_0(\n)\otimes\kk(\p)$ in terms
of the genus of $T$. These two bounds can be linked through
a comparison of the Euler-Poincar\'{e} characteristics of $T$
and its reduction, and it turns out that the number of $q^{2m}$-rational points on $X_0(\n)\otimes\kk(\p)$ is at least
$(q^m-1)(g-1)$, $g$ being the genus of $X_0(\n)\otimes\kk(\p)$. Thus,
it follows that the curves $X_0(\n)\otimes\kk(\p)$ are
assymptotically optimal when their genus increases.

\subsection{Two maps from $X_0(\n\p)$ to $X_0(\n)$}

Let $f_1:Y_0(\n\p)\to Y_0(\n)$ denote the forgetful morphism, the
map that ``forgets'' the level $\Gamma_0(\p)$-structure, that
associates with a triple $(\E,G,H)$ of an elliptic module $\E$,
a level $\Gamma_0(\n)$-structure $G$ and a level $\Gamma_0(\p)$-structure $H$, a pair $(\E,G)$. Denote by $f_2:Y_0(\n\p)\to Y_0(\n)$ the  morphism that associates with a triple $(\E,G,H)$ the pair
$(\E/H,GH/H)$ consisting of the quotient module (see \ref{quotientmodule}), together with its induced level
$\Gamma_0(\n)$-structure.

These are both maps between \emph{coarse} moduli spaces, so in order
to define them properly, one needs to pass to a fine cover first. For example, one can add a level $\Gamma(\m)$ structure, define the two
maps between $Y_0(\n\p)Y(\m)$ and $Y_0(\n)Y(\m)$, note that they
are $\GL(2,\A/\m)$-equivariant and pass to the quotient.
Both maps extend to maps $X_0(\n\p)\to X_0(\n)$ on the compactifications, since the generic fibres are smooth (\ref{allsmooth}). We will also denote the extensions by
$f_1$ and $f_2$.

Denote by $T$ the image of $X_0(\n\p)$ in $X_0(\n)\times X_0(\n)$
under $(f_1,f_2)$. Then $T$ is a two-dimensional integral scheme
over $\FF_q$, a closed subscheme of $X_0(\n)\times X_0(\n)$
and a flat scheme over $A$. Denote the normalisation
of $T$ by $\tilde{T}$. Then the generic fiber of $\tilde{T}\to\A$ is
necessarily smooth, since if it had a singularity, it would
extend to a singular curve on the surface $\tilde{T}$. 

\subsection{The reduction modulo $\p$}

Let $\Pi$ and $\Pi^t$ denote the graphs of the $q^m$-th power 
frobenius morphisms
 $\tau^m:X_0(\n)\otimes \kk(\p) \to X_0(\n)\otimes \kk(\p)$ and
 $\tau^m:X_0(\n)\otimes \kk(\p) \leftarrow X_0(\n)\otimes \kk(\p)$, respectively.  Consider $\Pi$, $\Pi^t$ and $\Pi\cup\Pi^t$ as reduced
closed subschemes of the surface $(X_0(\n)\otimes \kk(\p))\times(X_0(\n)\otimes \kk(\p))$ over $\kk(\p)$. 

\begin{proposition}
 $T\otimes_A \kk(\p)$  coincides with $\Pi\cup\Pi^t$.
\end{proposition}

\begin{proof}
 Let $\E$ be an elliptic module over an algebraic closure
 $\o{\kk(\p)}$ of
 $\kk(\p)$, $G$ a level $\Gamma_0(\n)$-structure on $\E$ and
 $H$ a level $\Gamma_0(\p)$-structure.
 Then by \ref{geomgammanought} there are two cases to consider.
 \begin{enumerate}
  \item $H$ is \emph{local}, then $H$ is the divisor  conisting
        of $q^m$ times the zero section. But this is precisely
        the kernel of the $q^m$-th power frobenius isogeny
        $\E\to \E'$, where $\E'$ is the elliptic
        module obtained
        by base extension $\tau^m:\o{\kk(\p)}\to \o{\kk(\p)}$. Hence
        $\E/H$ is isomorphic to $\E'$ and the geometric point
        corresponding to $(\E,G,H)$ maps under $(f_1,f_2)$ to the
        point corresponding to $((\E,G),(\E',G'))$, where $G'$
        is the level structure obtained by base extension of $G$.
  \item $H$ is \emph{\'{e}tale}. Denote by $\F$ the elliptic module
        $\E/H$. Let $\F'$ be the elliptic module obtained by
        base extension over $\tau^m$ as above. Then, as before,
        $\F'$ is isomorphic to the module $\F/q^m[0]$. This, in
        turn is isomorphic to $\E/\E[\p]$, and therefore also
        to $\E$. So the geometric point corresponding to $(\E,G,H)$
        maps under $(f_1,f_2)$ to the geometric point corresponding
        to $((\E,G),(\F,GH/H))$, or $((\F',(GH/H)'),(\F,GH/G))$.
   \end{enumerate}
  Thus, the image $T$ of $(f_1,f_2)$ coincides with
  $\Pi\cup\Pi^t$ on the open part
  $(Y_0(\n)\otimes \kk(\p))\times (Y_0(\n)\otimes \kk(\p))$ of 
  $(X_0(\n)\otimes \kk(\p))\times (X_0(\n)\otimes \kk(\p))$ and
  hence also on the entire surface.
\end{proof}

\begin{remark}
 The level $\Gamma_0(\n)$-structure in the above proposition
 can \emph{not} be replaced for a $\Gamma_1(\n)$-structure. To
 see this, take $H$ to be \'{e}tale in the above proof. Then
 $\E$ will still be isomorphic to $\F'$, and the induced
 level $\Gamma_1(\n)$-structure on $\F'$ will generate
 the same $\Gamma_0(\n)$-structure as before, but there is no
 way to decide which generator we get. 
\end{remark}

\subsection{Counting the rational points}

We are now in position to apply Ihara's trick in order to prove:

\begin{theorem}
 The number of $q^{2m}$-rational points of $X_0(\n)\otimes\kk(\p)$
 is at least $(q^m-1)(g-1)$, where $g$ equals the genus
 of $X_0(\n)\otimes \o{\kk(\p)}$.
\end{theorem}

Assuming the theorem, we can simply ``increase'' $\n$ to obtain
a tower of curves of increasing genus with assymptotically many
rational points over $\FF_{q^{2m}}$.

\begin{proof}
 (Almost identical to \cite[\S 1]{Ihara79}.)
 Let $g$ denote the genus of the generic fibre of $X_0(\n)$. Then,
 since $X_0(\n)\otimes \o{\kk(\p)}$ is smooth, $g$ is also the genus
 of the reduction. Denote by $g_0$ the genus of the generic
 fibre of $\tilde{T}$. Since both projections $T\to X_0(\n)$ have
 degree $q^m+1$, as can be seen on the fibre above $\p$, the Hurwitz
 formula gives
  $$ g_0 - 1 \geq (q^m+1)(g-1). $$
 
 The set of $\FF_{q^{2m}}$-rational points of
 $X_0(\n)\otimes\kk(\p)$ is in bijection with the set of geometric
 points in the intersection $\Pi \cap \Pi^t$, since it is the set
 of points invariant under $\tau^m\circ\tau^m$. 
 Such a point in the intersection $\Pi \cap \Pi^t$ is called
 \emph{special} if it is not a normal point of the two-dimensional
 scheme $T$. 
 
 The fibre above $\p$ of $\tilde{T}\to \A$ consists of
 two irreducible components, mapping to $\Pi$ and $\Pi^t$.
 These two components
 intersect above the special points of $\Pi\cap\Pi^t$. 
 Since the Euler-Poincar\'{e} characteristic of $\tilde{T}$ is
 constant along the fibers, we find 
  $$ g_0 - 1 = 2(g-1) + \#\{\text{special points}\}.$$
 Combined with the above inequality, this yields
  $$ \#\{\text{rational points over }\FF_{q^{2m}}\} \geq \
  \#\{\text{special points}\} \geq (q^m-1)(g-1). $$

\end{proof}

\bibliographystyle{plain}
\bibliography{../../master}

\def\cprime{$'$} \def\cprime{$'$} \def\cprime{$'$} \def\cprime{$'$}
  \def\cprime{$'$}
\begin{thebibliography}{10}

\bibitem{Deligne69}
P.~Deligne and D.~Mumford.
\newblock The irreducibility of the space of curves of given genus.
\newblock {\em Inst. Hautes \'Etudes Sci. Publ. Math.}, (36):75--109, 1969.

\bibitem{Dokchitser00}
Tim Dokchitser.
\newblock {\em Deformations of p-divisible groups}.
\newblock PhD thesis, University of Utrecht, 2000.

\bibitem{Drinfeld74E}
V.~G. Drinfel{\cprime}d.
\newblock Elliptic modules.
\newblock {\em Mat. Sb. (N.S.)}, 94(136):594--627, 656, 1974.

\bibitem{Elkies01}
Noam~D. Elkies.
\newblock Explicit towers of {D}rinfeld modular curves.
\newblock In {\em European Congress of Mathematics, Vol. II (Barcelona, 2000)},
  volume 202 of {\em Progr. Math.}, pages 189--198. Birkh\"auser, Basel, 2001.

\bibitem{EGA1}
A.~Grothendieck.
\newblock \'{E}l\'ements de g\'eom\'etrie alg\'ebrique. {I}. {L}e langage des
  sch\'emas.
\newblock {\em Inst. Hautes \'Etudes Sci. Publ. Math.}, (4):228, 1960.

\bibitem{EGA0}
A.~Grothendieck.
\newblock \'{E}l\'ements de g\'eom\'etrie alg\'ebrique. {0}.
\newblock {\em Inst. Hautes \'Etudes Sci. Publ. Math.}, (4):228,
  1960,1961,1963,1964,1965,1966,1967.

\bibitem{SGA1}
A.~Grothendieck.
\newblock {\em Rev\^etements \'etales et groupe fondamental ({SGA} 1)}.
\newblock Documents Math\'ematiques (Paris) [Mathematical Documents (Paris)],
  3. Soci\'et\'e Math\'ematique de France, Paris, 2003.
\newblock S\'eminaire de g\'eom\'etrie alg\'ebrique du Bois Marie 1960--61.
  [Algebraic Geometry Seminar of Bois Marie 1960-61], Directed by A.
  Grothendieck, With two papers by M. Raynaud, Updated and annotated reprint of
  the 1971 original [Lecture Notes in Math., 224, Springer, Berlin; MR0354651
  (50 \#7129)].

\bibitem{Grothendieck95}
Alexander Grothendieck.
\newblock Technique de descente et th\'eor\`emes d'existence en g\'eometrie
  alg\'ebrique. {I}. {G}\'en\'eralit\'es. {D}escente par morphismes
  fid\`element plats.
\newblock In {\em S\'eminaire Bourbaki, Vol.\ 5}, pages Exp.\ No.\ 190,
  299--327. Soc. Math. France, Paris, 1995.

\bibitem{Ihara79}
Yasutaka Ihara.
\newblock Congruence relations and {S}him\=ura curves.
\newblock In {\em Automorphic forms, representations and $L$-functions (Proc.
  Sympos. Pure Math., Oregon State Univ., Corvallis, Ore., 1977), Part 2},
  Proc. Sympos. Pure Math., XXXIII, pages 291--311. Amer. Math. Soc.,
  Providence, R.I., 1979.

\bibitem{Katz85}
Nicholas~M. Katz and Barry Mazur.
\newblock {\em Arithmetic moduli of elliptic curves}, volume 108 of {\em Annals
  of Mathematics Studies}.
\newblock Princeton University Press, Princeton, NJ, 1985.

\bibitem{Laumon96}
G{\'e}rard Laumon.
\newblock {\em Cohomology of {D}rinfeld modular varieties. {P}art {I}},
  volume~41 of {\em Cambridge Studies in Advanced Mathematics}.
\newblock Cambridge University Press, Cambridge, 1996.
\newblock Geometry, counting of points and local harmonic analysis.

\bibitem{Lehmkuhl00}
Thomas Lehmkuhl.
\newblock {\em Compactification of the Drinfeld modular surfaces
  (Habilitationsschrift)}.
\newblock G\"{o}ttingen, 2000.

\bibitem{Mumford94}
D.~Mumford, J.~Fogarty, and F.~Kirwan.
\newblock {\em Geometric invariant theory}, volume~34 of {\em Ergebnisse der
  Mathematik und ihrer Grenzgebiete (2) [Results in Mathematics and Related
  Areas (2)]}.
\newblock Springer-Verlag, Berlin, third edition, 1994.

\bibitem{Saidi96}
Mohamed Sa{\"{\i}}di.
\newblock Moduli schemes of {D}rinfeld modules.
\newblock In {\em Drinfeld modules, modular schemes and applications
  (Alden-Biesen, 1996)}, pages 17--31. World Sci. Publishing, River Edge, NJ,
  1997.

\bibitem{Tate97}
John Tate.
\newblock Finite flat group schemes.
\newblock In {\em Modular forms and Fermat's last theorem (Boston, MA, 1995)},
  pages 121--154. Springer, New York, 1997.

\bibitem{vanderHeiden03}
G.~J. van~der Heiden.
\newblock {\em Weil pairing and the Drinfeld modular curve}.
\newblock PhD thesis, University of Groningen, 2003.

\bibitem{Vladut91}
S.~G. Vl{\u{a}}du{\c{t}}.
\newblock {\em Kronecker's {J}ugendtraum and modular functions}, volume~2 of
  {\em Studies in the Development of Modern Mathematics}.
\newblock Gordon and Breach Science Publishers, New York, 1991.
\newblock Translated from the Russian by M. Tsfasman.

\bibitem{Vladut83}
S.~G. Vl{\`e}duts and V.~G. Drinfel{\cprime}d.
\newblock The number of points of an algebraic curve.
\newblock {\em Funktsional. Anal. i Prilozhen.}, 17(1):68--69, 1983.

\end{thebibliography}

\end{document}